\documentclass[11pt]{article}
\usepackage{amsfonts}
\usepackage{amsfonts}
\usepackage{amsfonts}
\usepackage{mathrsfs}
\usepackage{bm}
\usepackage{cite}
\usepackage[colorlinks,linkcolor=blue,citecolor=blue]{hyperref}
\usepackage{amssymb,amsmath} 
 \allowdisplaybreaks

\catcode`!=11
\let\!int\int \def\int{\displaystyle\!int}
\let\!lim\lim \def\lim{\displaystyle\!lim}
\let\!sum\sum \def\sum{\displaystyle\!sum}
\let\!sup\sup \def\sup{\displaystyle\!sup}
\let\!inf\inf \def\inf{\displaystyle\!inf}
\let\!cap\cap \def\cap{\displaystyle\!cap}
\let\!max\max \def\max{\displaystyle\!max}
\let\!min\min \def\min{\displaystyle\!min}
\let\!frac\frac \def\frac{\displaystyle\!frac}
\catcode`!=12

\let\oldsection\section
\renewcommand\section{\setcounter{equation}{0}\oldsection}

\allowdisplaybreaks
\def\pf{\it{Proof.}\rm\quad}

\def\N{\mathbb{N}}

\newtheorem{thm}{Theorem}[section]
\newtheorem{lem}[thm]{Lemma}
\newtheorem{cor}[thm]{Corollary}

\begin{document}
\title {\bf On harmonic numbers and nonlinear Euler sums}
\author{
{Ce Xu\thanks{Corresponding author. Email: xuce1242063253@163.com (C. XU)}\quad Yulin Cai\thanks{Yulin.Cai@math.u-bordeaux.fr (Y. Cai)}
}\\[1mm]
\small  $\ast$ School of Mathematical Sciences, Xiamen University\\
\small Xiamen
361005, P.R. China\\
\small  \dag Institut de Mathematiques de Bordeaux, Universite de Bordeaux\\
\small Bordeaux 33405, France
}

\date{}
\maketitle \noindent{\bf Abstract }
In this paper we are interested in Euler-type sums with products of harmonic numbers, Stirling numbers and Bell numbers. We discuss the analytic representations of Euler sums through values of polylogarithm function and Riemann zeta function. Moreover, we provide explicit evaluations for almost all Euler sums with weight $\leq 5$, which can be expressed in terms of zeta values and polylogarithms. Furthermore, we give explicit formula for several classes of Euler-related sums in terms of zeta values and harmonic numbers, and several examples are given. The given representations are new.
\\[2mm]
\\[2mm]
\tableofcontents
\section{Introduction and Preliminaries}

\subsection{Harmonic numbers and Euler sums}
The generalized $n$th-harmonic numbers of order $p$ are given by (\cite{FS1998})
\begin{align}\label{a1}
H_n^{(p)}:=\sum\limits_{j=1}^n\frac {1}{j^p},\quad  p \in \N,
\end{align}
where the $H_n\equiv H_n^{(1)}$ is the classical harmonic number (\cite{A1995}) and the empty sum $H_0^{(p)}$ is conventionally understood to be zero.

In response to a letter from Goldbach in 1742, Euler considered sums of the form (see Berndt \cite{B1985})
\[S_{p,q}:=\sum\limits_{n = 1}^\infty  {\frac{{H_n^{\left( p \right)}}}
{{{n^q}}}}, \]
where $p,q$ are positive integers with $q \geq 2$, and $w:=p+q$ denotes the weight of sums $S_{p,q}$. These sums are called the linear Euler sums today. Euler discovered that in the cases  $p=1,\ p=q$ and $p+q$ is less than 7 or when $p+q$ is odd and less than 13, the linear sums have evaluations in terms of zeta values, i.e., the values of the Riemann zeta function (\cite{A2000})
\[\zeta(s):=\sum\limits_{n = 1}^\infty {\frac {1}{n^{s}}},\Re(s)>1\]
at the positive integer arguments. Moreover, he conjectured that the double linear sums would be reducible to zeta values when $p+q$ is odd, and even gave what he hoped to obtain the general formula. In 1995, Borweins and Girgensohn \cite{BBG1995} proved conjecture and formula, and in 1994, Bailey, Borwein and Girgensohn \cite{BBG1994} conjectured that the linear sums $S_{p,q}$, when $p + q> 7$ and $p + q$ is even, are not reducible. Hence, the linear sums $S_{p,q}$ can be evaluated in terms of zeta values in the following cases: $p=1,p=q,p+q$ odd and $p+q=4,6$ with $q\geq 2$.

Similarly, the nonlinear Euler sums are the infinite sums whose general term is a product of harmonic numbers of index $n$ and a power of $n^{-1}$.
Let $\pi  := \left( {{\pi _1}, \ldots ,{\pi _k}} \right)$ be a partition of integer $p$ and $p = {\pi _1} +  \cdots  + {\pi _k}$ with ${\pi _1} \le {\pi _2} \le  \cdots  \le {\pi _k}$. The classical nonlinear Euler sum of index $\pi,q$ is defined as follows (see \cite{FS1998})
\begin{align}\label{a2}
{S_{\pi_1\pi_2\cdots\pi_k ,q}}: = \sum\limits_{n = 1}^\infty  {\frac{{H_n^{\left( {{\pi _1}} \right)}H_n^{\left( {{\pi _2}} \right)} \cdots H_n^{\left( {{\pi _k}} \right)}}}{{{n^q}}}},\quad q\geq 2
\end{align}
where the quantity $w:={\pi _1} +  \cdots  + {\pi _k} + q$ is called the weight and the quantity $k$ is called the degree. As usual, repeated summands in partitions are indicated by powers, so that for instance
\[{S_{{1^2}{2^3}4,q}} = {S_{112224,q}} = \sum\limits_{n = 1}^\infty  {\frac{{H_n^2[H^{(2)} _n]^3{H^{(4)} _n}}}{{{n^q}}}}. \]
The nonlinear Euler sums, i.e., $S_{\pi ,q}$ with $\pi$ having two or more parts, are more complicated. Such sums were already considered in \cite{BBG1994,FS1998,W2017,X2016,X2017,Xu2017,XC2017,Xu-2017}.
In \cite{FS1998}, Flajolet and Salvy gave an algorithm for reducing $S_{\pi_1\pi_2,q}$ to linear Euler sums when $\pi_1+\pi_2+q$ is even and $\pi_1,\pi_2,q>1$  (see Theorem 4.2 in the reference \cite{FS1998}. In \cite{X2016,X2017}, the first author jointly with Li, Yan and Shi proved that all quadratic Euler sums of the form
\[S_{1(p-1),q}=\sum\limits_{n = 1}^\infty  {\frac{{{H_n}{H_n^{(p-1)}}}}{{{n^q}}}}\ \ \left( {p,q\geq 2,\ p + q \le 10} \right)\]
can be reduced to polynomials in zeta values and linear sums. Moreover, in the recent papers \cite{Xu2017,W2017}, the first author \cite{Xu2017} shown that all Euler sums of weight $\leq 8$ are reducible to $\mathbb{Q}$-linear combinations of single zeta monomials with the addition of $\{S_{2,6}\}$ for weight 8. For weight 9, all Euler sums of the form ${S_{{\pi_1} \cdots {\pi_k},q}}$ with $q\in \{4,5,6,7\}$ are expressible polynomially in terms of zeta values. Wang and Lyu \cite{W2017} shown that all Euler sums of weight eight are reducible to linear sums, and proved that all Euler sums of weight nine are reducible to zeta values. In the most recent papers \cite{XC2017,Xu-2017,WC2017}, the first author has carried out fruitful cooperation with Wang proved that all Euler sums of weight ten can be expressed as a rational linear combination of $\zeta \left( {10} \right),{\zeta ^2}\left( 5 \right),\zeta \left( 3 \right)\zeta \left( 7 \right),\zeta \left( 2 \right)\zeta \left( 3 \right)\zeta \left( 5 \right),{\zeta ^2}\left( 3 \right)\zeta \left( 4 \right),\zeta \left( 2 \right){S_{2,6}}$ and ${S_{2,8}}$. All Euler sums of weight 11 are reducible to $S_{2,6},\ S_{18,2}$ and zeta values.

Similarly, in \cite{FS1998}, Flajolet and Salvy introduced the following alternating harmonic numbers
\begin{align}\label{a3}
\bar H_n^{\left( p \right)}: = \sum\limits_{j = 1}^n {\frac{{{{\left( { - 1} \right)}^{j - 1}}}}{{{j^p}}}} \quad{\rm and}\quad {{\bar H}_n}: = \bar H_n^{\left( 1 \right)},\quad p\in\N.
\end{align}
Clearly, we have the following relations for harmonic number and alternating harmonic number
\begin{align*}
  \bar H_n^{\left( p \right)}= H_{2\left[ {\left( {n + 1} \right)/2} \right] - 1}^{(p)} - \frac{1}
{{{2^p}}}\left( {H_{\left[ {n/2} \right]}^{(p)} + H_{\left[ {\left( {n - 1} \right)/2} \right]}^{(p)}} \right),\quad p\geq 2,
\end{align*}
\[{\bar H}_n = {H_n} - {H_{\left[ {n/2} \right]}},\]
where $[x]$ denotes the greatest integer less than or equal to $x$.\\
 The alternating harmonic number converges to the alternating Riemann zeta value ${\bar \zeta}(k)$:
 \[\mathop {\lim }\limits_{n \to \infty } {\bar H}_n^{\left( k \right)} = {\bar \zeta}( k),\quad k\in \N,\]
where the alternating Riemann zeta function is defined by (for more details, see \cite{FS1998})
$${\bar \zeta}(s):=\sum\limits_{n = 1}^\infty {\frac {(-1)^{n-1}}{n^{s}}},\quad \Re(s)\geq 1.$$
The generalized alternating Euler sums involving harmonic numbers and alternating harmonic numbers are defined by
\begin{align}
&{S_{\prod\limits_{i = 1}^l {{p_i}} \prod\limits_{j = 1}^k {{{\bar r}_j}} ,q}}: = \sum\limits_{n = 1}^\infty  {\frac{{H_n^{\left( {{p_1}} \right)} \cdots H_n^{\left( {{p_l}} \right)}\bar H_n^{\left( {{r_1}} \right)} \cdots \bar H_n^{\left( {{r_k}} \right)}}}{{{n^q}}}} ,\quad q>1 \label{a4}\\
&{S_{\prod\limits_{i = 1}^l {{p_i}} \prod\limits_{j = 1}^k {{{\bar r}_j}} ,\bar q}}: = \sum\limits_{n = 1}^\infty  {\frac{{H_n^{\left( {{p_1}} \right)} \cdots H_n^{\left( {{p_l}} \right)}\bar H_n^{\left( {{r_1}} \right)} \cdots \bar H_n^{\left( {{r_k}} \right)}}}{{{n^q}}}{{\left( { - 1} \right)}^{n - 1}}} ,\quad q\geq 1 \label{a5}
\end{align}
where the $\left( {{p_1}, \ldots ,{p_l}} \right)\in (\N)^l,\ \left( {{r_1}, \ldots ,{r_k}} \right)\in (\N)^k$ with $p_1\leq \cdots \leq p_l$ and $r_1\leq \cdots \leq r_k$ and the quantity $w: = {p_1} +  \cdots  + {p_l} + {r_1} +  \cdots  + {r_k} + q$ being called the weight and the quantity $l+k$ being the degree, $l\geq 0 $ and $ k\geq 1$. Since repeated summands in partitions are indicated by powers, we denote, for example, the sums
\[{S_{{1^3}2{{\bar 1}^2}\bar 3,5}} = \sum\limits_{n = 1}^\infty  {\frac{{H_n^3H_n^{\left( 2 \right)}\bar H_n^2\bar H_n^{\left( 3 \right)}}}{{{n^5}}}},\quad {S_{{1^2}3{{\bar 1}^3}\bar 3,\bar 4}} = \sum\limits_{n = 1}^\infty  {\frac{{H_n^2H_n^{\left( 3 \right)}\bar H_n^3\bar H_n^{\left( 3 \right)}}}{{{n^4}}}{{\left( { - 1} \right)}^{n - 1}}} .\]
It has been discovered in the course of the years that many alternating Euler sums admit expressions involving finitely the zeta values, $\ln2$ and polylogarithms. For example, in \cite{FS1998}, Flajolet and Salvy gave explicit reductions to zeta values and logarithm for all linear sums
\[{S_{p,q}} = \sum\limits_{n = 1}^\infty  {\frac{{H_n^{\left( p \right)}}}{{{n^q}}}} ,\:{S_{\bar p,q}} = \sum\limits_{n = 1}^\infty  {\frac{{\bar H_n^{\left( p \right)}}}{{{n^q}}}} ,\:{S_{p,\bar q}} = \sum\limits_{n = 1}^\infty  {\frac{{H_n^{\left( p \right)}}}{{{n^q}}}} {\left( { - 1} \right)^{n - 1}}\ {\rm and}\ {S_{\bar p,\bar q}} = \sum\limits_{n = 1}^\infty  {\frac{{\bar H_n^{\left( p \right)}}}{{{n^q}}}{{\left( { - 1} \right)}^{n - 1}}} \]
when $p+q$ is an odd weight. The evaluation of the above four linear sums in terms of values of the Riemann zeta function and polylogarithm function at positive integers is known when
 $(p,q) = (1,3), (2,2)$, or $p + q$ is odd (see \cite{BBG1995,C2010,FS1998,Xu2016}). Here are a few explicit evaluations
\begin{align*}
&S_{\bar 1, 3}= \frac{7}{4}\zeta \left( 3 \right)\ln 2 - \frac{5}{{16}}\zeta (4),\\
&S_{1,\bar 3}  =  - 2{\rm Li}{_4}\left( {1/2} \right) + \frac{{11{}}}{{4}} \zeta(4) + \frac{{{1}}}{{2}}\zeta(2){\ln ^2}2 - \frac{1}{{12}}{\ln ^4}2 - \frac{7}{4}\zeta \left( 3 \right)\ln 2,\\
&S_{\bar 1,\bar 3}{= \frac{3}{2}\zeta \left( 4 \right)}  + \frac{1}{2}\zeta \left( 2 \right){\ln ^2}2 - \frac{1}{{12}}{\ln ^4}2 - 2{\rm Li}{_4}\left( {1/2} \right),\\
&S_{2,\bar 2} =  - \frac{{51}}{{16}}\zeta (4) + 4{\rm{L}}{{\rm{i}}_4}\left( {1/2} \right) + \frac{7}{2}\zeta (3)\ln 2 - \zeta (2){\ln ^2}(2) + \frac{1}{6}{\ln ^4}2.
\end{align*}
Here ${\rm Li}_p(x)$ denotes the polylogarithm function, which is defined for $\left| x \right| \le 1$ by
\[{\rm Li}{_p}\left( x \right) := \sum\limits_{n = 1}^\infty  {\frac{{{x^n}}}{{{n^p}}}}, \Re(p)>1 .\]
Furthermore, we consider the nested sum
\[\sum\limits_{n = 1}^\infty  {\frac{{{y^n}}}{{{n^m}}}\left( {\sum\limits_{k = 1}^n {\frac{{{x^k}}}{{{k^p}}}} } \right)} ,\;x,y \in \left[ { - 1,1} \right),\; m,p \in \N.\]
By taking the sum over complementary pairs of summation indices, we obtain a simple reflection formula
\[\sum\limits_{n = 1}^\infty  {\frac{{{y^n}}}{{{n^m}}}\left( {\sum\limits_{k = 1}^n {\frac{{{x^k}}}{{{k^p}}}} } \right)}  + \sum\limits_{n = 1}^\infty  {\frac{{{x^n}}}{{{n^p}}}\left( {\sum\limits_{k = 1}^n {\frac{{{y^k}}}{{{k^m}}}} } \right)}  = {\rm Li}{_p}\left( x \right){\rm Li}{_m}\left( y \right) + {\rm Li}{_{p + m}}\left( {xy} \right).\]
Hence, we readily deduce the following four results of linear sums
\begin{align*}
&{S_{\bar 2,\bar 2} = } \frac{{13}}{{16}}\zeta \left( 4 \right),\quad S_{3,\bar 1}= \frac{{19}}{{16}}\zeta (4) - \frac{3}{4}\zeta \left( 3 \right)\ln 2,\\
& {S_{\bar 2, 2} = \frac{{85}}{{16}}\zeta \left( 4 \right)}  - 4{\rm{L}}{{\rm{i}}_4}\left( {1/2} \right) + \zeta \left( 2 \right){\ln ^2}2 - \frac{1}{6}{\ln ^4}2 - \frac{7}{2}\zeta \left( 3 \right)\ln 2,\\
&{S_{\bar 3,\bar 1} = } 2{\rm{L}}{{\rm{i}}_4}\left( {\frac{1}{2}} \right) + \frac{3}{4}\zeta \left( 3 \right)\ln 2 + \frac{1}{{12}}{\ln ^4}2 - \frac{1}{2}\zeta \left( 4 \right) - \frac{1}{2}\zeta \left( 2 \right){\ln ^2}2.
\end{align*}
However, in many other cases we are not able to obtain a formula for the Euler sum constant explicitly in terms of values of the Riemann zeta, logarithm and polylogarithm functions,
but we are able to obtain relations involving two of more Euler sum constants of the same weight.
In \cite{Xyz2016}, the first author jointly with Yang and Zhang proved that the alternating quadratic Euler sums
\begin{align*}
&{S_{{1^2},{\overline {2m}}}}: = \sum\limits_{n = 1}^\infty  {\frac{{H_n^2}}{{{n^{2m}}}}} {\left( { - 1} \right)^{n - 1}},\ {S_{{{\bar 1}^2},{\overline {2m}}}} = \sum\limits_{n = 1}^\infty  {\frac{{{\bar H}_n^2}}{{{n^{2m}}}}} {\left( { - 1} \right)^{n - 1}},\ {S_{1\bar 1,\overline {2m}}} = \sum\limits_{n = 1}^\infty  {\frac{{{H_n}{{\bar H}_n}}}{{{n^{2m}}}}} {\left( { - 1} \right)^{n - 1}},\\
&{S_{1\bar 1,2m}} = \sum\limits_{n = 1}^\infty  {\frac{{{H_n}{{\bar H}_n}}}{{{n^{2m}}}}} ,\ {S_{{{\bar 1}^2},2m}} = \sum\limits_{n = 1}^\infty  {\frac{{{H}_n^2}}{{{n^{2m}}}}}\ {\rm and} \ \sum\limits_{n = 1}^\infty  {\frac{{H_n^2}}{{{{\left( {2n + 1} \right)}^{2m}}}}}
\end{align*}
are reducible to polynomials in zeta values and to linear sums, and are able to give explicit values for certain of these sums in terms of the Riemann zeta function and the polylogarithm function, here $m$ is a positive integer. In \cite{BBG1994}, Bailey, Borwein and Girgensohn considered sums
\[\sum\limits_{n = 1}^\infty  {\bar H_n^q{{\left( {n + 1} \right)}^{ - p}}} ,\;\sum\limits_{n = 1}^\infty  {\bar H_n^q{{\left( {n + 1} \right)}^{ - p}}{{\left( { - 1} \right)}^{n - 1}}} \;{\rm{and}}\;\sum\limits_{n = 1}^\infty  {H_n^q{{\left( {n + 1} \right)}^{ - p}}{{\left( { - 1} \right)}^{n - 1}}} ,\quad \left( {p + q \le 5} \right)\]
and obtained a number of experimental identities using the experimental method, where $p\geq 2$ and $q\geq 1$ are positive integers.

Euler sums may be studied through a profusion of methods: combinatorial, analytic and al-
gebraic. There are many other researches on Euler sums and Euler type sums. Some related results for Euler sums may be seen in the works of \cite{BZB2008,BG1996,CH2017,C2003,F2005,CM1994,M2014,SX2017,S2015,SS2011,CX2017,XMZ2016}.

\subsection{Multiple harmonic (star) numbers and multiple zeta (star) values}
For $m\in \N,\ {\bf S}:=(s_1,s_2,...,s_m)\in (\N)^m$, and a non-negative integer $n$, the multiple harmonic number (MHN for short) and the multiple harmonic star number (MHSN  for short)
are defined by, respectively
\begin{align}
{\zeta _n}\left( {{s_1},{s_2}, \ldots ,{s_m}} \right) := \sum\limits_{1 \le {k_m} <  \cdots  < {k_1} \le n} {\frac{1}{{k_1^{{s_1}} \cdots k_m^{{s_m}}}}},\label{a6}\\
{\zeta^\star_n}\left( {{s_1},{s_2}, \ldots ,{s_m}} \right) := \sum\limits_{1 \le {k_m} \le  \cdots  \le {k_1} \le n} {\frac{1}{{k_1^{{s_1}} \cdots k_m^{{s_m}}}}}.\label{a7}
\end{align}
The integers $m$ and $\omega : = \left| {\bf S} \right|: = \sum\nolimits_{i = 1}^m {{s_i}}$ are called the depth and the weight of a multiple harmonic number or multiple harmonic star number. By convention, we put ${\zeta _n}\left( {\bf S} \right)=0$ if $n<m$, and ${\zeta _n}\left(\emptyset \right)={\zeta^\star _n}\left(\emptyset \right)=1$. By ${\left\{ {{s_1},{s_2}, \ldots ,{s_j}} \right\}^m}$ we denote the sequence of depth $mj$ with $m$ repetitions of ${\left( {{s_1},{s_2}, \ldots ,{s_j}} \right)}.$\\
If $s_1 > 1$, then the limits as $n\rightarrow \infty$ exist, and are the well-known multiple zeta and multiple star-zeta values ${\zeta}\left( {{s_1},{s_2}, \ldots ,{s_m}} \right)$ and ${\zeta^\star}\left( {{s_1},{s_2}, \ldots ,{s_m}} \right)$ respectively. We put a bar on top of $s_j\ (j=1,\cdots k)$ if there is a sign $(-1)^{n_j}$ appearing in the denominator on the right. For example,
\begin{align}
&{\zeta _n}\left( {{{\bar s}_1},{s_2},s_3, \ldots ,s_{m-1},{{\bar s}_m}} \right) = \sum\limits_{1 \le {k_m} <  \cdots  < {k_1} \le n} {\frac{{{{\left( { - 1} \right)}^{{k_1} + {k_m}}}}}{{k_1^{{s_1}} \cdots k_m^{{s_m}}}}},\label{a8}\\
&{\zeta^\star _n}\left( {{s_1},{{\bar s}_2},s_3 \ldots ,s_{m-1},{{\bar s}_m}} \right) = \sum\limits_{1 \le {k_m} \le  \cdots  \le {k_1} \le n} {\frac{{{{\left( { - 1} \right)}^{{k_2} + {k_m}}}}}{{k_1^{{s_1}} \cdots k_m^{{s_m}}}}} .\label{a9}
\end{align}
The sums of types (\ref{a8}) and (\ref{a9}) are called the alternating harmonic numbers and alternating harmonic star numbers, respectively. For alternating MHNs or MHSNs, we have the similar limit cases. For example
\begin{align}
&\zeta \left( {{{\bar s}_1},{s_2}, \ldots ,{{\bar s}_m}} \right) := \mathop {\lim }\limits_{n \to \infty } {\zeta _n}\left( {{{\bar s}_1},{s_2}, \ldots ,{{\bar s}_m}} \right) = \sum\limits_{1 \le {k_m} <  \cdots  < {k_1}} {\frac{{{{\left( { - 1} \right)}^{{k_1} + {k_m}}}}}{{k_1^{{s_1}} \cdots k_m^{{s_m}}}}},\label{a10}\\
&{\zeta ^ \star }\left( {{s_1},{{\bar s}_2}, \ldots ,{{\bar s}_m}} \right): = \mathop {\lim }\limits_{n \to \infty } \zeta _n^ \star \left( {{s_1},{{\bar s}_2}, \ldots ,{{\bar s}_m}} \right) =\sum\limits_{1 \le {k_m} \le  \cdots  \le {k_1} \le n} {\frac{{{{\left( { - 1} \right)}^{{k_2} + {k_m}}}}}{{k_1^{{s_1}} \cdots k_m^{{s_m}}}}}.\label{a11}
\end{align}
The sums of types (\ref{a10}) and (\ref{a11}) (one of more the $s_j$ barred) are called the alternating multiple zeta values and alternating multiple zeta star values, respectively. \\
 If the first index is barred, then the limits as $n\rightarrow \infty$ exist for all positive-integer indices, e.g.,
\[\zeta \left( {{\rm{\bar 1}}} \right) = \mathop {\lim }\limits_{n \to \infty } {\zeta _n}\left( {{\rm{\bar 1}}} \right) =  - \ln 2.\]
Multiple zeta values are certainly interesting and important. The number $\zeta(\bar 6,\bar 2)$ appeared in the quantum field theory literature in 1986 \cite{Br1986}, well before the the phrase ¡°multiple zeta values¡± had been coined (Multiple zeta values were introduced by the Hoffman \cite{H1992} and Zagier \cite{DZ1994} in the 1990s). The alternating multiple zeta values have been intensely studied since the early 1990s, and extensive tables of them were compiled, first by the group at Lille \cite{Bo2002}, and later as part of the Multiple Zeta Value Data Mine \cite{BBV2010}.

In the past several bases have been considered for both the multiple zeta values and alternating multiple zeta values. The vector space of multiple zeta values can be constructed allowing basis elements, which contain besides the zeta values the index of which is a Lyndon word products of this type of zeta values of lower weight. One basis of this kind is (see \cite{BBV2010,M2000,KP2013})
\begin{align*}
& w=2,\quad \zeta(2),\\
&w=3,\quad \zeta(3),\\
&w=4,\quad \zeta (4),\\
&w=5,\quad \zeta(5),\ \zeta(2)\zeta(3),\\
&w=6,\quad \zeta(6),\ \zeta^2(3),\\
&w=7,\quad \zeta(7),\ \zeta(2)\zeta(5),\ \zeta(3)\zeta(4),\\
&w=8,\quad \zeta(8),\ \zeta(3)\zeta(5),\ \zeta(2)\zeta^2(3),\ \zeta(5,3),\\
&w=9,\quad \zeta(9),\ \zeta(2)\zeta(7),\ \zeta(3)\zeta(6),\ \zeta(4)\zeta(5),\ \zeta^3(3),\\
&w=10,\quad \zeta(10),\ \zeta(3)\zeta(7),\ \zeta^2(5),\ \zeta(2)\zeta(3)\zeta(5),\ \zeta(4)\zeta^2(3),\ \zeta(2)\zeta(5,3),\ \zeta(7,3), etc.
\end{align*}
Rational relations among multiple zeta values are known through weight (sum of
the indices) 22 and tabulated in the Multiple Zeta Value Data Mine (henceforth
MZVDM) \cite{BBV2010}. \\
A basis for the $\mathbb{Q}$-vector space spanned by the set of alternating multiple zeta values for $1\leq w\leq 6$ should be given
\begin{align*}
&w=1,\quad \ln2,\\
&w=2,\quad \zeta(2),\ \ln^22,\\
&w=3,\quad \zeta(3),\ \zeta(2)\ln2,\ \ln^32,\\
&w=4,\quad {\rm Li}_4(1/2),\ \zeta(4),\ \zeta(3)\ln2,\ \zeta(2)\ln^22,\ \ln^42,\\
&w=5,\quad {\rm Li}_5(1/2),\ \zeta(5),\ {\rm Li}_4(1/2)\ln2,\ \zeta(2)\zeta(3),\ \zeta(4)\ln2,\ \zeta(3)\ln^22,\ \zeta(2)\ln^32,\ \ln^52,\\
&w=6,\quad {\rm Li}_6(1/2),\ \zeta(6),\ \zeta(\bar 5,\bar 1),\ {\rm Li}_5(1/2)\ln2,\ \zeta(5)\ln2,\ {\rm Li}_4(1/2)\zeta(2),\ {\rm Li}_4(1/2)\ln^22,\ \zeta^2(3),\\
&\quad\quad\quad\quad \zeta(2)\zeta(3)\ln2,\ \zeta(4)\ln^22,\ \zeta(3)\ln^32,\ \zeta(2)\ln^42,\ \ln^62, etc.
\end{align*}
Rational relations among alternating multiple zeta values are also tabulated in the MZVDM (through weight 12 \cite{BBV2010}).\\
The relations between multiple zeta (star) values and the values of the Riemann zeta function and polylogarithm function have been studied by many authors. For explicit evaluations of the multiple zeta (star) values or related alternating sums, readers may consult \cite{BBBL1997,BBBL2001,BB2003,EW2012,H1997,HO2003,CM1994,P2009,SC2012,Xu2017,DZ1994,DZ2012,Z2015}. For example, Zagier \cite{DZ2012} proved that the multiple zeta star values ${\zeta ^ \star \left( {{{\left\{ 2 \right\}}^a},3,{{\left\{ 2 \right\}}^b}} \right)} $ and multiple zeta values ${\zeta \left( {{{\left\{ 2 \right\}}^a},3,{{\left\{ 2 \right\}}^b}} \right)} $, where $a,b\in \N_0:=\N\cup \{0\}=\{0,1,2,\ldots\}$, are reducible to polynomials in zeta values, and gave explicit formulae. Hessami Pilehrood et al. \cite{KP2013,KP2017} and Li \cite{L2012} provide three new proofs of Zagier's formula for ${\zeta ^ \star \left( {{{\left\{ 2 \right\}}^a},3,{{\left\{ 2 \right\}}^b}} \right)}$ based on a finite identity for partial sums of the zeta-star series and hypergeometric series computations, respectively.

\subsection{Relations between Euler sums and multiple zeta values}

From the definitions of harmonic numbers and multiple harmonic numbers, we can find that
\[H_n^{\left( p \right)} = {\zeta _n}\left( p \right)\quad{\rm{and}}\quad\bar H_n^{\left( p \right)} =  - {\zeta _n}\left( {\bar p} \right),\quad p\in\N.\]
According to the rules of the ``harmonic algebra" or ``stuffle product", it is obvious that the products of any number of harmonic numbers and alternating harmonic number can be expressed in terms of (alternating) multiple harmonic numbers. For example,
\begin{align*}
&\bar H_n^2 = \zeta _n^2\left( {{\rm{\bar 1}}} \right) = \zeta \left( 2 \right) + 2{\zeta _n}\left( {{\rm{\bar 1,\bar 1}}} \right),\\
&{H_n}H_n^{\left( 3 \right)} = {\zeta _n}\left( {\rm{1}} \right){\zeta _n}\left( {\rm{3}} \right) = {\zeta _n}\left( {{\rm{1,3}}} \right) + {\zeta _n}\left( {{\rm{3,1}}} \right) + {\zeta _n}\left( {\rm{4}} \right),\\
&{H_n}{{\bar H}_n} =  - {\zeta _n}\left( {\rm{1}} \right){\zeta _n}\left( {{\rm{\bar 1}}} \right) =  - {\zeta _n}\left( {{\rm{\bar 2}}} \right) - {\zeta _n}\left( {{\rm{1,\bar 1}}} \right) - {\zeta _n}\left( {{\rm{\bar 1,1}}} \right).
\end{align*}
Note that when different indices are combined, the combination of a barred and unbarred index is barred, while the combination of two barred indices is unbarred.\\
The series
\[\sum\limits_{n = 1}^\infty  {\frac{{{\zeta _n}\left( {{i_1},{i_2}, \cdots ,{i_k}} \right)}}{{{n^p}}}} ,\quad i_j\in\N\]
can be written as the sum of multiple zeta value
\[\zeta \left( {p,{i_1},{i_2}, \cdots ,{i_k}} \right) + \zeta \left( {p + {i_1},{i_2}, \cdots ,{i_k}} \right).\]
Series of the form
\[\sum\limits_{n = 1}^\infty  {\frac{{{\zeta _n}\left( {{i_1},{i_2}, \cdots ,{i_k}} \right)}}{{{n^p}}}} \quad {\rm or}\quad \sum\limits_{n = 1}^\infty  {\frac{{{\zeta _n}\left( {{i_1},{i_2}, \cdots ,{i_k}} \right)}}{{{n^p}}}{{\left( { - 1} \right)}^n}} \]
with one of more of the $i_j$ barred, can be written as a sum of alternating multiple zeta values.\\
It is clear that every (alternating) Euler sum of weight $w$ and degree $n$ is clearly a $\mathbb{Q}$-linear combination of (alternating) multiple zeta values or (alternating) multiple zeta star values of weight $w$ and depth less than or equal to $n+1$. In other words, (alternating) multiple zeta (star) values are ``atomic" quantities into which (alternating) Euler sums decompose. For example,
\[{S_{{{\bar 1}^2},\bar 2}} =  - \zeta \left( {\bar 2,2} \right) + 2\zeta \left( {\bar 2,{\rm{\bar 1}},{\rm{\bar 1}}} \right) - \zeta \left( {\bar 4} \right) - 2\zeta \left( {3,\bar 1} \right),\]
\[{S_{13,4}} = \zeta \left( {{\rm{4,1,3}}} \right) + \zeta \left( {4,4} \right) + \zeta \left( {{\rm{4,3,1}}} \right) + \zeta \left( {5,3} \right) + \zeta \left( {\rm{8}} \right) + \zeta \left( {7,1} \right),\]
\[{S_{1{\rm{\bar 1,2}}}} =  - \zeta \left( {2,\bar 2} \right) - \zeta \left( {2,1,\bar 1} \right) - \zeta \left( {2,\bar 1,1} \right) - \zeta \left( {\bar 4} \right) - \zeta \left( {3,\bar 1} \right) - \zeta \left( {\bar 3,1} \right).\]
Hence, based on the above, we obtain the following two conclusion:

{\bf 1.} All Euler sums of weight $\in \{3,4,5,6,7,9\}$ can be expressed as a rational linear combination of products of zeta values. For weight 8, all such sums are the sum of a polynomial in zeta values and a rational multiple of $\zeta(5,3)$. For weight 10, all such sums are reducible to zeta values, $\zeta(5,3)$ and $\zeta(7,3)$. (All explicit formulas can be found in \cite{FS1998,X2016,W2017,WC2017,X2017,Xu2017})

{\bf 2.} All alternating Euler sums of weight $\leq 5$ can be expressed in terms of zeta values, polylogarithms and $\ln2$. (Apart of explicit formulas can be found in \cite{BBG1994,BBG1995,M2014,Xu2016,Xyz2016})

In this paper, we will give almost all explicit evaluations of alternating Euler sums with weight less or equal to $5$.

The evaluation of Euler sums and multiple zeta values has been useful in various areas of theoretical physics, including in
support of Feynman diagram calculations and in resolving open questions on Feynman diagram contributions and relations among special functions \cite{C2005}, including the dilogarithm, Clausen function, and generalized hypergeometric function \cite{C2008}. An array of Euler sums and multiple zeta values is required in calculations
of high energy physics. These quantities appear for instance in developing the scattering theory of massless quantum electrodynamics \cite{Bl1999}. Harmonic number sums also occur in computer science in the efficiency analysis of algorithms \cite{BBG1994}.

\subsection{Purpose and plan}
The purpose of this paper is to evaluate (alternating) Euler sums and related sums involving (alternating) harmonic numbers.

The remainder of this paper is organized as follows.

In the second section we evaluate some Euler related sums which involve harmonic numbers (or alternating harmonic numbers), Stirling numbers and Bell numbers. Then we use the evaluations obtained to establish many identities involving two or more Euler sums of the same weight, which can be expressed in terms of Riemann zeta values and polylogarithms.

In the third section we use certain integral representations to evaluate several series with harmonic numbers and alternating harmonic numbers. Moreover, we establish many relations involving two or more Euler sums of the same weight. Using these relations obtained we give many closed form of alternating Euler sums with weight $\leq 5$. Furthermore,  we also consider the following type of Euler-type sums
\[{S_{{m_1}{m_2}...{m_p}}}(r,k;{x_1}, \ldots ,{x_p};x) = \sum\limits_{n = 1}^\infty  {\frac{{{\zeta _n}\left( {{m_1};{x_1}} \right) \cdots {\zeta _n}\left( {{m_r};{x_p}} \right)}}{{(n + r)\left( {n + k} \right)}}}x^n,\quad p, k,r \in \N,\]
 and provide some explicit evaluation in a closed form in terms of zeta values and harmonic numbers.

\section{Series involving Stirling numbers, harmonic numbers and Bell numbers}

In this section, we will establish some explicit relationships that involve Euler sums and zeta values.

\subsection{The Bell polynomials, Stirling numbers and Bell numbers}

The exponential complete Bell polynomials $Y_n$ are defined by (see section 3.3 in the \cite{L1974} and section 2.8 in the \cite{R2002}).
\begin{align}\label{b1}
\exp \left( {\sum\limits_{m \ge 1}^{} {{x_m}\frac{{{t^m}}}{{m!}}} } \right) = 1 + \sum\limits_{k \ge 1}^{} {{Y_k}\left( {{x_1},{x_2}, \cdots ,{x_k}} \right)\frac{{{t^k}}}{{k!}}}.
\end{align}
Then ${Y_0}\left(  \cdot  \right) = 1$ and
\begin{align}\label{b2}
{Y_k}\left( {{x_1},{x_2}, \cdots ,{x_k}} \right) = \sum\limits_{{c_1} + 2{c_2} +  \cdots  + k{c_k} = k} {\frac{{k!}}{{{c_1}!{c_2}! \cdots {c_k}!}}{{\left( {\frac{{{x_1}}}{{1!}}} \right)}^{{c_1}}}{{\left( {\frac{{{x_2}}}{{2!}}} \right)}^{{c_2}}} \cdots {{\left( {\frac{{{x_k}}}{{k!}}} \right)}^{{c_k}}}} .
\end{align}
Additionally, by Eq. (2.44) of \cite{R2002}, the polynomials $Y_k(\cdot)$ satisfy the recurrence relation
\begin{align}\label{b3}
{Y_k}\left( {{x_1},{x_2}, \cdots ,{x_k}} \right) = \sum\limits_{j = 0}^{k - 1} {\left( {\begin{array}{*{20}{c}}
   {k - 1}  \\
   j  \\
\end{array}} \right){x_{k - j}}{Y_j}\left( {{x_1},{x_2}, \cdots ,{x_j}} \right)},\quad k\in\N.
\end{align}

The (unsigned) Stirling number of the first kind $s(n,k)$ is defined by \cite{L1974}
\begin{align}\label{b4}
n!x\left( {1 + x} \right)\left( {1 + \frac{x}{2}} \right) \cdots \left( {1 + \frac{x}{n}} \right) = \sum\limits_{k = 0}^n {s(n+1,k+1){x^{k + 1}}},
\end{align}
with $s(n,k)=0$, if $n<k$ and $s(n,0)=s(0,n)=0,\ s(0,0)=1$, or equivalently, by the generating function:
\begin{align}\label{b5}
{\ln ^k}\left( {1 - x} \right) = {\left( { - 1} \right)^k}k!\sum\limits_{n = k}^\infty  {s(n,k)\frac{{{x^n}}}{{n!}}} ,\;x \in \left[ { - 1,1} \right).
\end{align}
Moreover, the (unsigned) Stirling numbers $s(n,k)$ of the first kind satisfy a recurrence relation in the form
\begin{align}\label{b6}
s\left( {n,k} \right) = s\left( {n - 1,k - 1} \right) + \left( {n - 1} \right)s\left( {n - 1,k} \right).
\end{align}
By the definition of $s(n,k)$, we see that we may rewrite (\ref{b4}) as
\begin{align}\label{b7}
\sum\limits_{k = 0}^n {s(n+1,k+1)}x^{k}&= n!\exp \left\{ {\sum\limits_{j = 1}^n {\ln \left( {1 + \frac{x}{j}} \right)} } \right\}\nonumber\\
& = n!\exp \left\{ {\sum\limits_{j = 1}^n {\sum\limits_{k = 1}^\infty  {{{\left( { - 1} \right)}^{k - 1}}\frac{{{x^k}}}{{k{j^k}}}} } } \right\}\nonumber\\
& = n!\exp \left\{ {\sum\limits_{k = 1}^\infty  {{{\left( { - 1} \right)}^{k - 1}}\frac{{{H^{(k)} _n}{x^k}}}{k}} } \right\}.
\end{align}
Letting ${x_k} = {\left( { - 1} \right)^{k - 1}}\left( {k - 1} \right)!H_n^{\left( k \right)}$ in (\ref{b1}) and combining (\ref{b7}), the following relation holds:
\begin{align}\label{b8}
s\left( {n + 1,k + 1} \right) = \frac{{n!}}{{k!}}{Y_k}\left( {{H_n}, - 1!H_n^{\left( 2 \right)}, \cdots ,{{\left( { - 1} \right)}^{k - 1}}\left( {k - 1} \right)!H_n^{\left( k \right)}} \right).
\end{align}
Therefore, we know that $s(n,k)$ is a rational linear combination of products of harmonic numbers.
Then by (\ref{b3}), we have the recurrence
\begin{align}\label{b9}
s\left( {n + 1,k + 1} \right) = \frac{1}{k}\sum\limits_{j = 0}^{k - 1} {{{\left( { - 1} \right)}^{k - j - 1}}H_n^{\left( {k - j} \right)}s\left( {n + 1,j + 1} \right)},\quad k\in\N,
\end{align}
which further yields
\begin{align*}
&s(n+1,1)= {n}!,\ s(n+1,2) ={n}!{H_{n }},\ s(n+1,3) = \frac{{n!}}{2}\left[ {H_{n}^2 - {H^{(2)}_{n}}} \right],\\
&s(n+1,4)= \frac{{n!}}{6}\left[ {H_{n}^3 - 3{H_{n}}{H^{(2)}_{n}} + 2{H^{(3)} _{n}}} \right], \\
&s(n+1,5) = \frac{{n!}}{{24}}\left[ {H_{n}^4 - 6{H^{(4)} _{n}} - 6H_{n}^2{H^{(2)} _{n}} + 3(H^{(2)} _{n})^2 + 8H_{n}^{}{H^{(3)} _{n}}} \right].
\end{align*}

Throughout this article we define a Bell number $\{Y_k(n)\}$ by the Bell polynomials
\begin{align}\label{b10}
{Y_k}\left( n \right) = {Y_k}\left( {{H_n},1!H_n^{\left( 2 \right)}, \cdots ,\left( {k - 1} \right)!H_n^{\left( k \right)}} \right),
\end{align}
Thus, according to the recurrence (\ref{b3}), we obtain
\[{Y_0}\left( n \right) = 1,\quad {Y_k}\left( n \right) = \sum\limits_{j = 0}^{k - 1} {\left( {\begin{array}{*{20}{c}}
   {k - 1}  \\
   j  \\
\end{array}} \right)\left( {k - j - 1} \right)!H_n^{\left( {k - j} \right)}{Y_j}\left( n \right)},\]
from which, we know that ${Y_k}\left( n \right)$ is a rational linear combination of products of harmonic numbers, and we have
$${Y_1}\left( n \right) = {H_n},{Y_2}\left( n \right) = H_n^2 + {H^{(2)} _n},{Y_3}\left( n \right) =  H_n^3+ 3{H_n}{H^{(2)} _n}+ 2{H^{(3)} _n},$$
\[{Y_4}\left( n \right) = H_n^4 + 8{H_n}{H^{(3)} _n} + 6H_n^2{H^{(2)} _n} + 3(H^{(2)} _n)^2 + 6{H^{(4)} _n},\]

\subsection{Some lemmas}
 The following lemmas will be useful in the development of the main theorems.
\begin{lem}(\cite{Xu2017})\label{lem1} For integers $m\geq 0$ and $n\geq1$, then the following integral identity holds:
\begin{align}\label{b11}
\int\limits_0^x {{t^{n - 1}}{{\left( {\ln t} \right)}^m}} dt = \sum\limits_{l = 0}^m {l!\left( {\begin{array}{*{20}{c}}
   m  \\
   l  \\
\end{array}} \right)\frac{{{{\left( { - 1} \right)}^l}}}{{{n^{l + 1}}}}{{\left( {\ln x} \right)}^{m - l}}{x^n}},\ x\in (0,1).
\end{align}
\end{lem}

\begin{lem}\label{lem2} For positive integers $m$ and nonnegative integers $k$,
\begin{align}\label{b12}
W\left( {m,k} \right) = \int\limits_0^1 {\frac{{{{\ln }^m}\left( {1 - x} \right){{\ln }^k}\left( x \right)}}{x}} dx,
\end{align}
then
\begin{align}\label{b13}
W\left( {m,k} \right) = {\left( { - 1} \right)^{m + k}}m!k!\zeta \left( {k + 2,{{\left\{ 1 \right\}}^{m - 1}}} \right).
\end{align}
\end{lem}
\pf  Expanding out $\ln^m(1-x)$ and using the lemma \ref{lem1},
\begin{align*}
W\left( {m,k} \right) &= {\left( { - 1} \right)^m}\sum\limits_{{i_1} = 1}^\infty  {\sum\limits_{{i_2} = 1}^\infty  { \cdots \sum\limits_{{i_m} = 1}^\infty  {\int\limits_0^1 {\frac{{{x^{{i_1} + {i_2} +  \cdots  + {i_m} - 1}}{{\ln }^k}\left( x \right)}}{{{i_1}{i_2} \cdots {i_m}}}dx} } } } \\
& = {\left( { - 1} \right)^m}\sum\limits_{{i_1} = 1}^\infty  {\sum\limits_{{i_2} = 1}^\infty  { \cdots \sum\limits_{{i_m} = 1}^\infty  {\frac{{{{\left( { - 1} \right)}^k}k!}}{{{i_1}{i_2} \cdots {i_m}{{\left( {{i_1} +  \cdots  + {i_m}} \right)}^{k + 1}}}}} } }.
\end{align*}
By Lemma 4.3 in the \cite{H1992},
\[\sum\limits_{{i_1} = 1}^\infty  {\sum\limits_{{i_2} = 1}^\infty  { \cdots \sum\limits_{{i_m} = 1}^\infty  {\frac{1}{{{i_1}{i_2} \cdots {i_m}{{\left( {{i_1} +  \cdots  + {i_m}} \right)}^{k + 1}}}}} } }  = \sum\limits_{{n_1} > {n_2} >  \cdots  > {n_m} \ge 1} {\frac{{m!}}{{n_1^{k + 2}{n_2} \cdots {n_k}}}} ,\]
and the conclusion follows.\hfill$\square$\\
In \cite{BBBL1997},  Borwein, Bradley and Broadhurst gave the generating function
\begin{align}\label{b14}
\sum\limits_{n,m = 1}^\infty  {\zeta \left( {m + 1,{{\left\{ 1 \right\}}^{n - 1}}} \right){x^m}{y^n} = 1 - \exp \left( {\sum\limits_{n = 2}^\infty  {\zeta \left( n \right)\frac{{{x^n} + {y^n} - {{\left( {x + y} \right)}^n}}}{n}} } \right)} ,
\end{align}
which implies that for any $m,\ n\in \N$, the multiple zeta value ${\zeta \left( {m + 1,{{\left\{ 1 \right\}}^{n - 1}}} \right)}$ can be
represented as a polynomial of zeta values with rational coefficients, and we have the duality formula
\[\zeta \left( {n + 1,{{\left\{ 1 \right\}}^{m - 1}}} \right) = \zeta \left( {m + 1,{{\left\{ 1 \right\}}^{n - 1}}} \right).\]
In particular, one can find explicit formulas for small weights.
\[\begin{array}{l}
 \zeta \left( {2,{{\left\{ 1 \right\}}^m}} \right) = \zeta \left( {m + 2} \right), \\
 \zeta \left( {3,{{\left\{ 1 \right\}}^m}} \right) = \frac{{m + 2}}{2}\zeta \left( {m + 3} \right) - \frac{1}{2}\sum\limits_{k = 1}^m {\zeta \left( {k + 1} \right)\zeta \left( {m + 2 - k} \right)} . \\
 \end{array}\]
In \cite{X2016}, the first author jointly with Yan and Shi gave the following recurrence relation of $W(m,k)$
\begin{align}\label{b15}
W\left( {m{\rm{,}}k} \right){\rm{ = }} - {\psi ^{\left( {m + k} \right)}}\left( 1 \right)/\left( {k + 1} \right) - \sum\limits_{i = 1}^{m - 1} {\sum\limits_{j = 0}^{k - 1} {\left( {\begin{array}{*{20}{c}}
   {m - 1}  \\
   i  \\
\end{array}} \right)} } \left( {\begin{array}{*{20}{c}}
   k  \\
   j  \\
\end{array}} \right)W\left( {i{\rm{,}}j} \right){\psi ^{\left( {m + k - i - j - 1} \right)}}\left( 1 \right),
\end{align}
where $\psi \left( z \right)$ denotes the digamma function (or called Psi function ) which is defined as the logarithmic derivative of the well known gamma function:
\[\psi \left( z \right) := \frac{d}{{dz}}\left( {\ln \Gamma \left( z \right)} \right) = \frac{{\Gamma '\left( z \right)}}{{\Gamma \left( z \right)}}.\]
From (\ref{b14}) or (\ref{b15}), we see that the values of integrals $W(m,k)$ can be expressed as a rational linear combination of products of zeta values.

\begin{lem}\label{lem3} For positive integers $m$ and nonnegative integers $k$,
\begin{align}\label{b16}
V\left( {m,k} \right) = \int\limits_0^1 {\frac{{{{\ln }^m}\left( {1 + x} \right){{\ln }^k}\left( x \right)}}{x}} dx,
\end{align}
then
\begin{align}\label{b17}
V\left( {m,k} \right) = {\left( { - 1} \right)^{m + k}}m!k!\zeta \left( {\overline{k + 2},{{\left\{ 1 \right\}}^{m - 1}}} \right).
\end{align}
\end{lem}
\pf The proof of Lemma \ref{lem3} is similar to the proof of Lemma \ref{lem2}.\hfill$\square$

\begin{lem}\label{lem4}
If $m\geq1$ is a integer and $z \in \left[ {0,1} \right]$
, then we have
\begin{align}\label{b18}
\int\limits_0^z {\frac{{{{\ln }^m}\left( {1 + x} \right)}}{x}} dx = &\frac{1}{{m + 1}}{\ln ^{m + 1}}\left( {1 + z} \right) + m!\left( {\zeta \left( {m + 1} \right) - {\rm Li}{_{m + 1}}\left( {\frac{1}{{1 + z}}} \right)} \right)\nonumber \\&- m!\sum\limits_{j = 1}^m {\frac{{{{\ln }^{m - j + 1}}\left( {1 + z} \right)}}{{\left( {m - j + 1} \right)!}}} {\rm Li}{_j}\left( {\frac{1}{{1 + z}}} \right).
\end{align}
\end{lem}
\pf We note that the integral in (\ref{b18}), which can be rewritten as
\begin{align}\label{b19}
\int\limits_0^z {\frac{{{{\ln }^m}\left( {1 + x} \right)}}{x}} dx\mathop  = \limits^{t = 1 + x}  & \int\limits_1^{1 + z} {\frac{{{{\ln }^m}t}}{{t - 1}}} dt\mathop  = \limits^{u = {t^{ - 1}}} {\left( { - 1} \right)^{m + 1}}\int\limits_1^{{{(1 + z)^{-1}}}} {\frac{{{{\ln }^m}u}}{{u - {u^2}}}} du\nonumber \\
           =& {\left( { - 1} \right)^{m + 1}}\left\{ {\int\limits_1^{{{(1 + z)^{-1}}}} {\frac{{{{\ln }^m}u}}{u}} du + \int\limits_1^{{{(1 + z)^{-1}}}} {\frac{{{{\ln }^m}u}}{{1 - u}}} du} \right\}\nonumber \\
           =&\frac{1}{{m + 1}}{\ln ^{m + 1}}\left( {1 + z} \right) + {\left( { - 1} \right)^{m + 1}}\int\limits_1^{{{(1 + z)^{-1}}}} {\frac{{{{\ln }^m}u}}{{1 - u}}} du .
\end{align}
Then with help of formula (\ref{b11}) we may deduce the desired result. \hfill$\square$

\begin{lem}(see \cite{X2016})
For integers $n\geq1$ and $k\geq0$, then
\begin{align}\label{b20}
\int\limits_0^1 {{t^{n - 1}}{{\ln }^k}\left( {1 - t} \right)} dt = {\left( { - 1} \right)^k}\frac{{{Y_k}\left( n \right)}}{n},
\end{align}
where the Bell number ${Y_k}\left( n \right)$ is defined by (\ref{b10}), and ${Y_k}\left( n \right) = k!\zeta _n^ \star \left( {{{\{ 1\} }^k}} \right)  $ .
\end{lem}

\begin{lem}(\cite{B1989,G2016})\label{lem6}  For positive integer $n,m$ and any suitable sequences $a_i\in \mathbb{C} \ (i=1,2\ldots,n)$,
\begin{align*}
&\prod\limits_{i = 1}^n\frac{1}{{ {\left( {1 - {a_i}t} \right)} }} = \sum\limits_{m = 0}^\infty  {{A_m}\left( n \right){t^m}} ,\quad {A_0}\left( n \right) = 1,\\
&\prod\limits_{i = 1}^n {\left( {1 + {a_i}t} \right)}  = \sum\limits_{m = 0}^\infty  {{{\bar A}_m}\left( n \right){t^m}},\quad {{\bar A}_0}\left( n \right) = 1.
\end{align*}
then
\begin{align}
{A_m}\left( n \right): = \sum\limits_{1 \le {k_m} \le  \cdots  \le {k_1} \le n} {{a_{{k_1}}} \cdots {a_{{k_m}}}} ,\label{b21}\\
{{\bar A}_m}\left( n \right): = \sum\limits_{1 \le {k_m} <  \cdots  < {k_1} \le n} {{a_{{k_1}}} \cdots {a_{{k_m}}}} .\label{b22}
\end{align}
\end{lem}
Here
$\left| t \right| < \min \left\{ {{{\left| {{a_1}} \right|}^{ - 1}}, \cdots ,{{\left| {{a_n}} \right|}^{ - 1}}} \right\}.$

\subsection{Main results and proofs}\label{sec2.3}

\begin{thm}\label{thm2.7} For positive integer $m$, then the following symmetric-function identities hold:
\begin{align}\label{b23}
{h_m} = \frac{1}{m}\sum\limits_{i = 0}^{m - 1} {{h_i}{p_{m - i}}} \quad {\rm and}\quad {e_m} = \frac{{{{\left( { - 1} \right)}^{m - 1}}}}{m}\sum\limits_{i = 0}^{m - 1} {{{\left( { - 1} \right)}^i}{e_i}{p_{m - i}}} ,
\end{align}
where $h_m=h_m(x_1,x_2,\cdots,x_n)$ and $e_m=e_m(x_1,x_2,\cdots,x_n)$ are respectively the complete and elementary symmetric functions of degree $m$, and $p_m=x_1^m+x_2^m+\cdots+x_n^m$ is the $m$th power sum (for details introductions, see \cite{H2017}, the notation follows \cite{Ma1995}).
\end{thm}
\pf From Hoffman's paper \cite{H2017}, we know that
\[E\left( t \right) = \sum\limits_{j = 0}^\infty  {{e_j}{t^j}}  = \prod\limits_{i = 1}^n {\left( {1 + t{x_i}} \right)} \quad {\rm and}\quad H\left( t \right) = \sum\limits_{j = 0}^\infty  {{h_j}{t^j}}  = \prod\limits_{i = 1}^n {\frac{1}{{1 - t{x_i}}}} \]
are the respective generating functions of the elementary and complete symmetric functions.
Note that we can rewrite $E(t)$ and $H(t)$ as
\[E\left( t \right) = \exp \left( {\sum\limits_{k = 1}^\infty  {\frac{{{{\left( { - 1} \right)}^{k - 1}}}}{k}{p_k}{t^k}} } \right)\quad {\rm and}\quad H\left( t \right) = \exp \left( {\sum\limits_{k = 1}^\infty  {\frac{1}{k}{p_k}{t^k}} } \right).\]
Here $\left| t \right| < \min \left\{ {{{\left| {{x_1}} \right|}^{ - 1}}, \cdots ,{{\left| {{x_n}} \right|}^{ - 1}}} \right\}.$
Letting $x_k=(-1)^{k-1}(k-1)!p_k$ and $(k-1)!p_k$ in (\ref{b1}), then comparing the coefficients of $t^k$, we can find that
\[{e_k} = \frac{{{Y_k}\left( {{p_1}, - 1!{p_2}, \cdots ,{{\left( { - 1} \right)}^{k - 1}}\left( {k - 1} \right)!{p_k}} \right)}}{{k!}}\quad{\rm and}\quad {h_k} = \frac{{{Y_k}\left( {{p_1},1!{p_2}, \cdots ,\left( {k - 1} \right)!{p_k}} \right)}}{{k!}}.\]
Using the recurrence relation (\ref{a3}), by a simple calculation, we obtain the desired results.\hfill$\square$\\
(Theorem \ref{thm2.7} can also be found in Entry 12(a) of Chapter 14 in Ramanujan's second notebook \cite{B1989}).
Hence, the elementary symmetric function $e_k$ and the complete symmetric function $h_k$ can be expressed in terms of the power sums, i.e., (see \cite{H2017})
\[{e_k} = {P_k}\left( {{p_1},{p_2}, \cdots ,{p_k}} \right)\quad {\rm and}\quad {h_k} = Q\left( {{p_1},{p_2}, \cdots ,{p_k}} \right).\]
Moreover, from the generating function $E(t)$ and $H(t)$, we have
\begin{align*}
 \sum\limits_{m = 0}^\infty  {{e_m}\left( {{x_1}, \cdots ,{x_n}} \right){t^m}}  &= \prod\limits_{i = 1}^n {\left( {1 + t{x_i}} \right)}  = \frac{1}{{1 + t{x_{n + 1}}}}\sum\limits_{m = 0}^\infty  {{e_m}\left( {{x_1}, \cdots ,{x_n},{x_{n + 1}}} \right){t^m}}  \\
  &= \sum\limits_{m = 0}^\infty  {\left( {\sum\limits_{k = 0}^m {{{\left( { - 1} \right)}^k}x_{n + 1}^k{e_{m - k}}\left( {{x_1}, \cdots ,{x_{n + 1}}} \right)} } \right){t^m}},  \\
 \sum\limits_{m = 0}^\infty  {{h_m}\left( {{x_1}, \cdots ,{x_n},{x_{n + 1}}} \right){t^m}}  &= \prod\limits_{i = 1}^{n + 1} {\frac{1}{{1 - t{x_i}}}}  = \frac{1}{{1 - t{x_{n + 1}}}}\sum\limits_{m = 0}^\infty  {{h_m}\left( {{x_1}, \cdots ,{x_n},} \right){t^m}}  \\
  &= \sum\limits_{m = 0}^\infty  {\left( {\sum\limits_{k = 0}^m {x_{n + 1}^k{h_{m - k}}\left( {{x_1}, \cdots ,{x_n}} \right)} } \right){t^m}}.
\end{align*}
Then by considering the coefficient of $t^m$, we can get the following relations
\begin{align*}
 {e_m}\left( {{x_1}, \cdots ,{x_n}} \right) = \sum\limits_{k = 0}^m {{{\left( { - 1} \right)}^k}x_{n + 1}^k{e_{m - k}}\left( {{x_1}, \cdots ,{x_{n + 1}}} \right)} , \\
 {h_m}\left( {{x_1}, \cdots ,{x_{n + 1}}} \right) = \sum\limits_{k = 0}^m {x_{n + 1}^k{h_{m - k}}\left( {{x_1}, \cdots ,{x_n}} \right)} .
\end{align*}
The polynomials $P_n$ and $Q_n$ can be written as determinants (see \cite{HO2003})
\[n!{P_n}\left( {{y_1}, \ldots ,{y_n}} \right) = \left| {\begin{array}{*{20}{c}}
   {{y_1}} & 1 & 0 &  \cdots  & 0  \\
   {{y_2}} & {{y_1}} & 2 &  \cdots  & 0  \\
    \vdots  &  \vdots  &  \vdots  &  \ddots  &  \vdots   \\
   {{y_{n - 1}}} & {{y_{n - 2}}} & {{y_{n - 3}}} &  \cdots  & {n - 1}  \\
   {{y_n}} & {{y_{n - 1}}} & {{y_{n - 2}}} &  \cdots  & {{y_1}}  \\
\end{array}} \right|\]
and
\[n!{Q_n}\left( {{y_1}, \ldots ,{y_n}} \right) = \left| {\begin{array}{*{20}{c}}
   {{y_1}} & { - 1} & 0 &  \cdots  & 0  \\
   {{y_2}} & {{y_1}} & { - 2} &  \cdots  & 0  \\
    \vdots  &  \vdots  &  \vdots  &  \ddots  &  \vdots   \\
   {{y_{n - 1}}} & {{y_{n - 2}}} & {{y_{n - 3}}} &  \cdots  & { - n + 1}  \\
   {{y_n}} & {{y_{n - 1}}} & {{y_{n - 2}}} &  \cdots  & {{y_1}}  \\
\end{array}} \right|.\]
Here are a few explicit evaluations; some of them were given previously by Hoffman and Ohno.
\begin{align*}
& {P_1}\left( {{y_1}} \right) = {y_1},\quad {P_2}\left( {{y_1},{y_2}} \right) = \frac{{y_1^2 - {y_2}}}{2}, \\
& {P_3}\left( {{y_1},{y_2},{y_3}} \right) = \frac{{y_1^3 - 3{y_1}{y_2} + 2{y_3}}}{6}, \\
& {P_4}\left( {{y_1},{y_2},{y_3},{y_4}} \right) = \frac{{y_1^4 - 6y_1^2{y_2} + 8{y_1}{y_3} + 3y_2^2 - 6{y_4}}}{{24}}, \\
& {P_5}\left( {{y_1},{y_2},{y_3},{y_4},{y_5}} \right) = \frac{{y_1^5 - 10y_1^3{y_2} + 20y_1^2{y_3} + 15{y_1}y_2^2 - 30{y_1}{y_4} - 20{y_2}{y_3} + 24{y_5}}}{{120}} \\
\end{align*}
and
\begin{align*}
& {Q_1}\left( {{y_1}} \right) = {y_1},\quad {Q_2}\left( {{y_1},{y_2}} \right) = \frac{{y_1^2 + {y_2}}}{2}, \\
& {Q_3}\left( {{y_1},{y_2},{y_3}} \right) = \frac{{y_1^3 + 3{y_1}{y_2} + 2{y_3}}}{6}, \\
 &{Q_4}\left( {{y_1},{y_2},{y_3},{y_4}} \right) = \frac{{y_1^4 + 6y_1^2{y_2} + 8{y_1}{y_3} + 3y_2^2 + 6{y_4}}}{{24}}, \\
 &{Q_5}\left( {{y_1},{y_2},{y_3},{y_4},{y_5}} \right) = \frac{{y_1^5 + 10y_1^3{y_2} + 20y_1^2{y_3} + 15{y_1}y_2^2 + 30{y_1}{y_4} + 20{y_2}{y_3} + 24{y_5}}}{{120}}. \\
\end{align*}
By Lemma {\ref{lem6}}, we have the relations
\begin{align}
&{P_m}\left( {{p_1}, \ldots ,{p_m}} \right) = \underset{k_1 = 1}{\overset{n}{\sum}}x_{k_1}\underset{k_2 = 1}{\overset{k_1-1}{\sum}}x_{k_2}\cdots\underset{k_m = 1}{\overset{k_{m-1}-1}{\sum}}x_{k_m},\label{b24}\\
&{Q_m}\left( {{p_1}, \ldots ,{p_m}} \right) =  \underset{k_1 = 1}{\overset{n}{\sum}}x_{k_1}\underset{k_2 = 1}{\overset{k_1}{\sum}}x_{k_2}\cdots\underset{k_m = 1}{\overset{k_{m-1}}{\sum}}x_{k_m}.\label{b25}
\end{align}
Hence, according to the definitions of Stirling numbers of the first kind $s(n,k)$ and Bell numbers $Y_k(n)$, we can get the well known identities
\begin{align}
& s\left( {n,k} \right) = {P_{k - 1}}\left( {{H_{n - 1}}, \ldots ,H_{n - 1}^{\left( {k - 1} \right)}} \right) = \left( {n - 1} \right)!{\zeta _{n - 1}}\left( {{{\left\{ 1 \right\}}^{k - 1}}} \right),\label{b26}\\
& {Y_k}\left( n \right) = k!{Q_k}\left( {{H_n}, \ldots ,H_n^{\left( k \right)}} \right) = k!\zeta _n^ \star \left( {{{\{ 1\} }^k}} \right),\label{b27}
\end{align}
which were proved in \cite{Xu2017}.\\
Therefore, by (\ref{b26}), then
\begin{align}\label{b28}
\sum\limits_{n = k}^\infty  {\frac{{s\left( {n,k} \right)}}{{n!{n^p}}}}  = \sum\limits_{n = p}^\infty  {\frac{{s\left( {n,p} \right)}}{{n!{n^k}}}}  = \zeta \left( {p + 1,{{\left\{ 1 \right\}}^{k - 1}}} \right).
\end{align}
The two series are equal, which is an instance of duality for multiple zeta values. (See also the discussion on pp. 526-7 of Chen \cite{CH2017}).

Using Theorem \ref{thm2.7} and formulas (\ref{b24}) and (\ref{b25}), we are able to obtain many relations involving multiple zeta-star values and zeta values. For example, taking ${x_i} = \frac{1}{{{i^p}}},p>1\;\left( {i = 1, \cdots ,n} \right)$ in Theorem \ref{thm2.7}, (\ref{b24}) and (\ref{b25}), then letting $n \to \infty $, we deduce the results
\begin{align*}
&\zeta \left( {{{\left\{ p \right\}}^m}} \right) = \frac{{{{\left( { - 1} \right)}^{m - 1}}}}{m}\sum\limits_{i = 0}^{m - 1} (-1)^i{\zeta \left( {{{\left\{ p \right\}}^i}} \right)\zeta \left( {pm - pi} \right)},\\
&\zeta^\star \left( {{{\left\{ p \right\}}^m}} \right) = \frac{1}{m}\sum\limits_{i = 0}^{m - 1} {\zeta^\star \left( {{{\left\{ p \right\}}^i}} \right)\zeta \left( {pm - pi} \right)}.
\end{align*}
Setting $x_i=\frac {(-1)^i}{i}$ and $m=6$ in (\ref{b24}), then letting $n \to \infty $ yields
\begin{align*}
 \zeta \left( {{{\left\{ {\bar 1} \right\}}^6}} \right) = & - \frac{5}{{128}}\zeta \left( 6 \right) + \frac{1}{{32}}{\zeta ^2}\left( 3 \right) + \frac{3}{{16}}\zeta \left( 5 \right)\ln 2 - \frac{1}{8}\zeta \left( 2 \right)\zeta \left( 3 \right)\ln 2 + \frac{1}{{32}}\zeta \left( 4 \right){\ln ^2}2 \\
  &+ \frac{1}{{24}}\zeta \left( 3 \right){\ln ^3}2 - \frac{1}{{48}}\zeta \left( 2 \right){\ln ^4}2 + \frac{1}{{720}}{\ln ^6}2.
\end{align*}

Now we establish some relations between zeta values and series involving Stirling numbers, harmonic numbers and Bell numbers.

\begin{thm}(\cite{Xu2017}) \label{thm2.8}
For positive integers $k$ and $p$, we have
\begin{align}
& \sum\limits_{n = k}^\infty  {\frac{{s\left( {n,k} \right){Y_p}\left( n \right)}}{{n!n}}}  = p!\left( {\begin{array}{*{20}{c}}
   {k + p}  \\
   p  \\
\end{array}} \right)\zeta \left( {k + p + 1} \right),\label{b29}\\
& \sum\limits_{n = k}^\infty  {\frac{{s\left( {n,k} \right)}}{{n!n}}{H^{(p+1)} _n}}  = \sum\limits_{n = p}^\infty  {\frac{{s\left( {n,p} \right)}}{{n!n}}H^{(k+1)}_n} ,\label{b30}\\
& p!\sum\limits_{n = p}^\infty  {\frac{{s\left( {n,p} \right){Y_k}\left( n \right)}}{{n!{n^2}}}}  + k!\sum\limits_{n =k}^\infty  {\frac{{s\left( {n,k} \right){Y_p}\left( n \right)}}{{n!{n^2}}}}  = k!p!\zeta \left( {k + 1} \right)\zeta \left( {p + 1} \right).\label{b31}
 \end{align}
\end{thm}
\pf Multiplying (\ref{b5}) by $x^{-1}\ln^p(1-x)$ and integrating over (0,1), using (\ref{b20}), we deduce that
\[\int\limits_0^1 {\frac{{{{\ln }^{p + k}}\left( {1 - x} \right)}}{x}} dx{\rm{ }} = {\left( { - 1} \right)^{k + p}}k!\sum\limits_{n = k}^\infty  {\frac{{s\left( {n,k} \right){Y_p}\left( n \right)}}{{n!n}}} \]
By Lemma \ref{lem2}, the integral on the left hand side is equal to
\[\int\limits_0^1 {\frac{{{{\ln }^{p + k}}\left( {1 - x} \right)}}{x}} dx = {\left( { - 1} \right)^{p + k}}\left( {p + k} \right)!\zeta \left( {p + k + 1} \right).\]
Summing these two contributions yields the result (\ref{b29}).

To prove (\ref{b30}), we need the well known identity
\begin{align}\label{b32}
\int\limits_0^1 {\frac{{{x^n}{{\ln }^p}x}}{{1 - x}}} dx = {\left( { - 1} \right)^p}p!\left( {\zeta \left( {p+ 1} \right) - H_n^{\left( {p + 1} \right)}} \right),\quad p\in\N.
\end{align}
Replacing $x$ by $t$ in (\ref{b5}), then dividing (\ref{b5}) by $t$ and integrating over the interval $(0,x)$, the result is
\begin{align}\label{b33}
\int\limits_0^x {\frac{{{{\ln }^k}\left( {1 - t} \right)}}{t}} dt = {\left( { - 1} \right)^k}k!\sum\limits_{n = k}^\infty  {\frac{{s\left( {n,k} \right)}}{{n!n}}{x^n}} ,\quad - 1 \le x \le 1.
\end{align}
Multiplying (\ref{b33}) by $\frac{{{{\ln }^p}x}}{{1 - x}}$ and integrating over the interval $(0,1)$, by (\ref{b32}), we obtain
\begin{align}\label{b34}
\int\limits_0^1 {\frac{{{{\ln }^p}x}}{{1 - x}}} \int\limits_0^x {\frac{{{{\ln }^k}\left( {1 - t} \right)}}{t}} dtdx = {\left( { - 1} \right)^{k + p}}k!p!\sum\limits_{n = k}^\infty  {\frac{{s\left( {n,k} \right)}}{{n!n}}\left( {\zeta \left( {p + 1} \right) - H_n^{\left( {p + 1} \right)}} \right)} .
\end{align}
Noting that the integral on the left hand side of this equation can be written as
\[\int\limits_0^1 {\frac{{{{\ln }^p}x}}{{1 - x}}} \int\limits_0^x {\frac{{{{\ln }^k}\left( {1 - t} \right)}}{t}} dtdx\mathop  = \limits_{t = 1 - y}^{x = 1 - u} \int\limits_0^1 {\frac{{{{\ln }^p}\left( {1 - u} \right)}}{u}} \int\limits_u^1 {\frac{{{{\ln }^k}y}}{{1 - y}}} dydu = \int\limits_0^1 {\frac{{{{\ln }^p}\left( {1 - x} \right)}}{x}} \int\limits_x^1 {\frac{{{{\ln }^k}t}}{{1 - t}}} dtdx.\]
By integration by parts, we see that
\begin{align*}
 \int\limits_0^1 {\frac{{{{\ln }^p}\left( {1 - x} \right)}}{x}} \int\limits_x^1 {\frac{{{{\ln }^k}t}}{{1 - t}}} dtdx&= \int\limits_0^1 {\left( {\int\limits_x^1 {\frac{{{{\ln }^k}t}}{{1 - t}}} dt} \right)} d\left( {\int\limits_0^x {\frac{{{{\ln }^p}\left( {1 - t} \right)}}{t}dt} } \right) \\
  &= \left. {\left( {\int\limits_x^1 {\frac{{{{\ln }^k}t}}{{1 - t}}} dt} \right)\left( {\int\limits_0^x {\frac{{{{\ln }^p}\left( {1 - t} \right)}}{t}dt} } \right)} \right|_0^1 + \int\limits_0^1 {\frac{{{{\ln }^k}x}}{{1 - x}}} \int\limits_0^x {\frac{{{{\ln }^p}\left( {1 - t} \right)}}{t}} dtdx \\
  &= {\left( { - 1} \right)^{p + k}}p!k!\sum\limits_{n = k}^\infty  {\frac{{s\left( {n,p} \right)}}{{n!n}}\left( {\zeta \left( {k + 1} \right) - {H^{(k+1)}_n}} \right)} .
 \end{align*}
Then with the help of formula (\ref{b28}) we easily obtain the result.

To prove identity (\ref{b31}), we consider the following integral
\begin{align*}
 &\int\limits_0^1 {\int\limits_0^x {\frac{{{{\ln }^k}\left( {1 - x} \right){{\ln }^p}\left( {1 - t} \right)}}{{xt}}} dtdx}  = {\left( { - 1} \right)^p}p!\sum\limits_{n = p}^\infty  {\frac{{s\left( {n,p} \right)}}{{nn!}}\int\limits_0^1 {{x^{n - 1}}{{\ln }^k}\left( {1 - x} \right)} dx}  \\
  &= {\left( { - 1} \right)^{p + k}}p!\sum\limits_{n = p}^\infty  {\frac{{s\left( {n,p} \right){Y_k}\left( n \right)}}{{{n^2}n!}}}  \\
  &= \left. {\left( {\int\limits_0^x {\frac{{{{\ln }^p}\left( {1 - t} \right)}}{t}dt} } \right)\left( {\int\limits_0^x {\frac{{{{\ln }^k}\left( {1 - t} \right)}}{t}dt} } \right)} \right|_0^1 - \int\limits_0^1 {\int\limits_0^x {\frac{{{{\ln }^k}\left( {1 - t} \right){{\ln }^p}\left( {1 - x} \right)}}{{tx}}} dtdx}  \\
 & = {\left( { - 1} \right)^{k + p}}k!p!\zeta \left( {k+ 1} \right)\zeta \left( {p + 1} \right) - {\left( { - 1} \right)^{k + p}}k!\sum\limits_{n = k}^\infty  {\frac{{s\left( {n,k} \right){Y_p}\left( n \right)}}{{{n^2}n!}}}.\
 \end{align*}
We thus arrive at (\ref{b31}) to complete the proof of Theorem \ref{thm2.8}.\hfill$\square$

\begin{thm}\label{thm2.9} For positive integer $m$, then
\begin{align}\label{b35}
\sum\limits_{n = m}^\infty  {\frac{{{{\bar H}_n}s\left( {n,m} \right)}}{{nn!}}}  = \zeta \left( {2,{{\left\{ 1 \right\}}_{m - 1}}} \right)\ln 2 - \zeta \left( {\bar 2,{{\left\{ 1 \right\}}_{m - 1}}} \right)\ln 2 - \left( {m + 1} \right)\zeta \left( {\bar 2,{{\left\{ 1 \right\}}_{m }}} \right).
\end{align}
\end{thm}
\pf Using integration by parts we deduce
\begin{align}\label{b36}
\int_0^x {{t^{n - 1}}\ln \left( {1 - t} \right)} dt = \frac {1}{n}\left\{{x^n}\ln \left( {1 - x} \right) - \sum\limits_{j = 1}^n {\frac{{{x^j}}}{j}}  - \ln \left( {1 - x} \right)\right\},\quad  -1\leq x<1.
\end{align}
Taking $x=-1$ in (\ref{b36}) yields
\[\int\limits_0^1 {{x^{n - 1}}\ln \left( {1 + x} \right)} dx = \frac{{{{\left( { - 1} \right)}^n}}}{n}\left\{ {{{\left( { - 1} \right)}^n}\ln 2 - \ln 2 + {{\bar H}_n}} \right\}.\]
Next, we consider the following integral
\begin{align*}
 \int\limits_0^1 {\frac{{{{\ln }^{m + 1}}\left( {1 + x} \right)}}{x}} dx &= {\left( { - 1} \right)^m}m!\sum\limits_{n = m}^\infty  {s\left( {n,m} \right)\frac{{{{\left( { - 1} \right)}^n}}}{{n!}}} \int\limits_0^1 {{x^{n - 1}}\ln \left( {1 + x} \right)} dx \\
  &= {\left( { - 1} \right)^m}m!\ln 2\sum\limits_{n = m}^\infty  {\frac{{s\left( {n,m} \right)}}{{nn!}}{{\left( { - 1} \right)}^n}}  - {\left( { - 1} \right)^m}m!\ln 2\sum\limits_{n = m}^\infty  {\frac{{s\left( {n,m} \right)}}{{nn!}}}  \\&\quad + {\left( { - 1} \right)^m}m!\sum\limits_{n = m}^\infty  {\frac{{{{\bar H}_n}s\left( {n,m} \right)}}{{nn!}}}
\end{align*}
By Lemma \ref{lem2} and Lemma \ref{lem3}
\begin{align*}
&\int\limits_0^1 {\frac{{{{\ln }^m}\left( {1 + x} \right)}}{x}} dx = {\left( { - 1} \right)^m}m!\sum\limits_{n = m}^\infty  {\frac{{s\left( {n,m} \right)}}{{nn!}}{{\left( { - 1} \right)}^n}}=(-1)^mm!\zeta\left(\bar 2,\{1\}_{m-1}\right),\\
&\int\limits_0^1 {\frac{{{{\ln }^m}\left( {1 - x} \right)}}{x}} dx = {\left( { - 1} \right)^m}m!\sum\limits_{n = m}^\infty  {\frac{{s\left( {n,m} \right)}}{{nn!}}}=(-1)^mm!\zeta\left( 2,\{1\}_{m-1}\right),
\end{align*}
Thus, combining the above identities, we complete the proof of Theorem \ref{thm2.9}.\hfill$\square$\\
By Lemma {\ref{lem3}} and Lemma \ref{lem4}, we can get the well known alternating multiple zeta values (\cite{BBBL1997})
\begin{align*}
\zeta \left( {\bar 2,{{\left\{ 1 \right\}}_{m-1}}} \right) =& \frac{{{{\left( { - 1} \right)}^{m }}}}{{\left( {m + 1} \right)!}}{\ln ^{m + 1}}2 + {\left( { - 1} \right)^{m }}\left( {\zeta \left( {m + 1} \right) - {\rm{L}}{{\rm{i}}_{m + 1}}\left( {\frac{1}{2}} \right)} \right) \\&- {\left( { - 1} \right)^{m }}\sum\limits_{j = 1}^{m } {\frac{{{{\left( {\ln 2} \right)}^{m + 1 - j}}}}{{\left( {m + 1- j} \right)!}}} {\rm{L}}{{\rm{i}}_j}\left( {\frac{1}{2}} \right).
\end{align*}

\begin{thm}
For positive integer $m$, then
\begin{align}
&\sum\limits_{n = 1}^\infty  {\frac{{{Y_m}\left( n \right)}}{n}{{\left( { - 1} \right)}^{n - 1}}}  = m!{\rm{L}}{{\rm{i}}_{m + 1}}\left( {\frac{1}{2}} \right),\label{b37}\\
&\sum\limits_{n = 1}^\infty  {\frac{{{Y_{m - 1}}(n){{\bar H}_n}}}{n}{{\left( { - 1} \right)}^{n - 1}}}  = m!\left( {\zeta \left( {m + 1} \right) - {\rm{L}}{{\rm{i}}_{m + 1}}\left( {\frac{1}{2}} \right)} \right),\label{b38}\\
&\sum\limits_{n = m}^\infty  {\frac{{{{\bar H}_n}s\left( {n,m} \right)}}{{n!}}{{\left( { - 1} \right)}^{n - 1}}}  = {\rm{L}}{{\rm{i}}_{m + 1}}\left( {\frac{1}{2}} \right) - \frac{{{{\left( { - 1} \right)}^m}}}{{m!}}{\ln ^{m + 1}}2.\label{b39}
\end{align}
\end{thm}
\pf To prove the first and second identities. Using (\ref{b20}), by a direct calculation, it is easy to see that
\begin{align}\label{b40}
\int\limits_0^1 {\frac{{{{\ln }^m}\left( {1 - x} \right)}}{{1 + x}}} dx =& \sum\limits_{n = 1}^\infty  {{{\left( { - 1} \right)}^{n - 1}}} \int\limits_0^1 {{x^{n - 1}}{{\ln }^m}\left( {1 - x} \right)} dx = {\left( { - 1} \right)^m}\sum\limits_{n = 1}^\infty  {\frac{{{Y_m}\left( n \right)}}{n}{{\left( { - 1} \right)}^{n - 1}}}.
\end{align}
Applying the change of variable $t\rightarrow 1-x$ to integral above and using (\ref{b11}), we obtain
\[\int\limits_0^1 {\frac{{{{\ln }^m}\left( {1 - t} \right)}}{{1 + t}}} dt = {\left( { - 1} \right)^m}m!{\rm{L}}{{\rm{i}}_{m + 1}}\left( {\frac{1}{2}} \right).\]
Thus, we deduce (\ref{b37}). Using the Cauchy product formula of power series, we can find that
\begin{align}\label{b41}
\frac{{\ln \left( {1 - x} \right)}}{{1 + x}} = \sum\limits_{n = 1}^\infty  (-1)^n{{{\bar H}_n}{{{ x} }^n}} ,\quad - 1 < x < 1.
\end{align}
Therefore, taking into account formulas (\ref{b20}) and (\ref{b41}), we obtain
\begin{align*}
 \int\limits_0^1 {\frac{{{{\ln }^m}\left( {1 - x} \right)}}{{x\left( {1 + x} \right)}}} dx &= \int\limits_0^1 {\frac{{{{\ln }^m}\left( {1 - x} \right)}}{x}} dx - \int\limits_0^1 {\frac{{{{\ln }^m}\left( {1 - x} \right)}}{{1 + x}}} dx \\
 & = \sum\limits_{n = 1}^\infty  {{{\bar H}_n}{{\left( { - 1} \right)}^n}\int\limits_0^1 {{x^{n - 1}}{{\ln }^{m - 1}}\left( {1 - x} \right)} dx}  \\
  &= {\left( { - 1} \right)^m}\sum\limits_{n = 1}^\infty  {\frac{{{Y_{m - 1}}(n){{\bar H}_n}}}{n}{{\left( { - 1} \right)}^{n - 1}}}.
 \end{align*}
Then with the help of formulas (\ref{b13}) and (\ref{b40}) we may deduce the result  (\ref{b38}).

To prove the third identity. By using (\ref{b5}), we have
\begin{align}\label{b42}
 \int\limits_0^1 {\frac{{{{\ln }^m}\left( {1 - x} \right)}}{{1 + x}}} dx &= {\left( { - 1} \right)^m}m!\sum\limits_{n = m}^\infty  {\frac{{s\left( {n,m} \right)}}{{n!}}} \int\limits_0^1 {\frac{{{x^n}}}{{1 + x}}} dx\nonumber \\
  &= {\left( { - 1} \right)^m}m!\sum\limits_{n = m}^\infty  {\frac{{{{\bar H}_n}s\left( {n,m} \right)}}{{n!}}{{\left( { - 1} \right)}^{n - 1}}}  + {\ln ^{m + 1}}2 .
 \end{align}
Comparing (\ref{b40}) with (\ref{b42}) completes the proof of (\ref{b39}).\hfill$\square$

\subsection{Some identities on Euler sums}

\begin{align}
 &\sum\limits_{n = 1}^\infty  {\frac{{{H_n}}}{n}{{\left( { - 1} \right)}^{n - 1}}}  = \frac{{\zeta \left( 2 \right) - {{\ln }^2}2}}{2},\label{b43} \\
&\sum\limits_{n = 1}^\infty  {\frac{{H_n^2 + {H^{(2)} _n}}}{n}{{\left( { - 1} \right)}^{n - 1}}}  = \frac{7}{8}\zeta \left( 3 \right) - \frac{1}{2}\zeta \left( 2 \right)\ln 2 + \frac{1}{6}{\ln ^3}2,\label{b44}\\
&\sum\limits_{n = 1}^\infty  {\frac{{{H_n}{{\bar H}_n}}}{{{n^2}}}}  = \frac{{43}}{{16}}\zeta (4) + \frac{3}{4}\zeta (2){\ln ^2}2 - \frac{1}{8}{\ln ^4}2 - 3{\rm{L}}{{\rm{i}}_4}\left( {\frac{1}{2}} \right),\label{b45}\\
& \sum\limits_{n = 1}^\infty  {\frac{{H_n^3 + 3{H_n}{H^{(2)} _n} + 2{H^{(3)} _n}}}{n}{{\left( { - 1} \right)}^{n - 1}}}  = 6{\rm{L}}{{\rm{i}}_4}\left( {\frac{1}{2}} \right),\label{b46}\\
&\sum\limits_{n = 1}^\infty  {\frac{{\bar H_n^3 + 3{{\bar H}_n}H_n^{\left( 2 \right)} + 2\bar H_n^{\left( 3 \right)}}}{n}{{\left( { - 1} \right)}^{n - 1}}}  = \frac{1}{4}{\ln ^4}2 + \frac{3}{2}\zeta \left( 3 \right)\ln 2 + \frac{{27}}{8}\zeta \left( 4 \right) + \frac{3}{2}\zeta \left( 2 \right){\ln ^2}2,\label{b47}\\
 &\sum\limits_{n = 1}^\infty  {\frac{{\left( {H_n^2 + {H^{(2)} _n}} \right){{\bar H}_n}}}{n}{{\left( { - 1} \right)}^{n - 1}}}  = 6\left( {\zeta \left( 4 \right) - {\rm{L}}{{\rm{i}}_4}\left( {\frac{1}{2}} \right)} \right),\label{b48}\\
&\sum\limits_{n = 1}^\infty  {\frac{{\left( {H_n^2 - H_n^{\left( 2 \right)}} \right){{\bar H}_n}}}{n}{{\left( { - 1} \right)}^{n - 1}}}  - 2\sum\limits_{n = 1}^\infty  {\frac{{{H_n}{{\bar H}_n}}}{{{n^2}}}{{\left( { - 1} \right)}^{n - 1}}}  + 2\sum\limits_{n = 1}^\infty  {\frac{{{{\bar H}_n}}}{{{n^3}}}{{\left( { - 1} \right)}^{n - 1}}} \nonumber\\& = 2{\rm{L}}{{\rm{i}}_4}\left( {\frac{1}{2}} \right) + \frac{1}{3}{\ln ^4}2,\label{b49}\\
&\sum\limits_{n = 1}^\infty  {\frac{{\left( {H_n^3 + 3{H_n}H_n^{\left( 2 \right)} + 2H_n^{\left( 3 \right)}} \right){{\bar H}_n}}}{n}{{\left( { - 1} \right)}^{n - 1}}}  = 24\left( {\zeta \left( 5 \right) - {\rm{L}}{{\rm{i}}_5}\left( {\frac{1}{2}} \right)} \right),\label{b50}\\
&\sum\limits_{n = 1}^\infty  {\frac{{H_n^4 + 8{H_n}H_n^{\left( 3 \right)} + 6H_n^2H_n^{\left( 2 \right)} + 3{{\left( {H_n^{\left( 2 \right)}} \right)}^2} + 6H_n^{\left( 4 \right)}}}{n}{{\left( { - 1} \right)}^{n - 1}}}  = 24{\rm{L}}{{\rm{i}}_5}\left( {\frac{1}{2}} \right),\label{b51}\\
&3\sum\limits_{n = 1}^\infty  {\frac{{\left( {H_n^2 - H_n^{\left( 2 \right)}} \right){{\bar H}_n}}}{{{n^2}}}}  - 6\sum\limits_{n = 1}^\infty  {\frac{{{H_n}{{\bar H}_n}}}{{{n^3}}}}  + 6\sum\limits_{n = 1}^\infty  {\frac{{{{\bar H}_n}}}{{{n^4}}}} \nonumber \\
& = 12\zeta \left( 4 \right)\ln 2 - \frac{5}{2}\zeta \left( 2 \right){\ln ^3}2 + \frac{{11}}{{20}}{\ln ^5}2 + \frac{{21}}{4}\zeta \left( 3 \right){\ln ^2}2 + 18{\rm{L}}{{\rm{i}}_4}\left( {\frac{1}{2}} \right)\ln 2 \nonumber \\ &\quad+ 24{\rm{L}}{{\rm{i}}_5}\left( {\frac{1}{2}} \right) - 24\zeta \left( 5 \right),\label{b52}\\
&\sum\limits_{n = 1}^\infty  {\frac{{{\bar H}_n^4 + 8{{\bar H}_n}\bar H_n^{\left( 3 \right)} + 6\bar H_n^2H_n^{\left( 2 \right)} + 3{{\left( {H_n^{\left( 2 \right)}} \right)}^2} + 6H_n^{\left( 4 \right)}}}{n}{{\left( { - 1} \right)}^{n - 1}}} \nonumber\\
& = \frac{9}{2}\zeta \left( 5 \right) + 3\zeta \left( 2 \right)\zeta \left( 3 \right) + \frac{{27}}{2}\zeta \left( 4 \right)\ln 2 + 3\zeta \left( 3 \right){\ln ^2}2 + 2\zeta \left( 2 \right){\ln ^3}2 + \frac{1}{5}{\ln ^5}2,\label{b53}\\
&\sum\limits_{n = 1}^\infty  {\frac{{\left( {H_n^3 - 3{H_n}H_n^{\left( 2 \right)} + H_n^{\left( 3 \right)}} \right){{\bar H}_n}}}{n}{{\left( { - 1} \right)}^{n - 1}}}  - 3\sum\limits_{n = 1}^\infty  {\frac{{\left( {H_n^2 - H_n^{\left( 2 \right)}} \right){{\bar H}_n}}}{{{n^2}}}{{\left( { - 1} \right)}^{n - 1}}} \nonumber\\
&\quad + 6\sum\limits_{n = 1}^\infty  {\frac{{{H_n}{{\bar H}_n}}}{{{n^3}}}{{\left( { - 1} \right)}^{n - 1}}}  - 6\sum\limits_{n = 1}^\infty  {\frac{{{{\bar H}_n}}}{{{n^4}}}{{\left( { - 1} \right)}^{n - 1}}}  = 6{\rm{L}}{{\rm{i}}_5}\left( {\frac{1}{2}} \right) - \frac{1}{4}{\ln ^5}2.\label{b54}
\end{align}
These identities can be obtained from the main theorems which are presented in the section \ref{sec2.3}. For example, setting $m=1,2,3,4$ in (\ref{b37}) respectively yields (\ref{b43}), (\ref{b44}), (\ref{b46}) and (\ref{b51}). Taking ${x_k} = \frac{{{{\left( { - 1} \right)}^{k - 1}}}}{k}$ and $m=4,5$ in (\ref{b25}) and Theorem \ref{thm2.7}, then letting $n\rightarrow \infty$, we obtain the results (\ref{b47}) and (\ref{b53}). etc

\section{Some explicit evaluation of Euler sums}
In the section, by constructing the integrals of polylogarithm functions, we establish some relations involving two or more Euler sums of the same weight. Then we use the relations obtained to evaluate alternating Euler sums, and give explicit formulas. First, we need the following lemma.

\begin{lem}(\cite{Xyz2016})\label{lem3.1} For positive integers $n,m$ and real $x\in [-1,1)$, we have
\begin{align}\label{c1}
\int\limits_0^x {{t^{n - 1}}{\rm Li}{_m}\left( t \right)dt}  =& \sum\limits_{i = 1}^{m - 1} {{{\left( { - 1} \right)}^{i - 1}}\frac{{{x^n}}}{{{n^i}}}{\rm Li}{_{m + 1 - i}}\left( x \right)} - \frac{{{{\left( { - 1} \right)}^m}}}{{{n^m}}}\left( {\sum\limits_{k = 1}^n {\frac{{{x^k}}}{k}} } \right)\nonumber\\
& + \frac{{{{\left( { - 1} \right)}^m}}}{{{n^m}}}\ln \left( {1 - x} \right)\left( {{x^n} - 1} \right).
\end{align}
\end{lem}

\subsection{Main theorems and corollaries}\label{sec3.1}
\begin{thm}\label{thm3.2}
Let $m,p$ be integers with $m,p>0$,  $\left| x \right|,\left| y \right|,\left| z \right| \in \left( {0,1} \right)$, we have
\begin{align}\label{c2}
 &\sum\limits_{n = 1}^\infty  {\frac{{{\rm{L}}{{\rm{i}}_m}\left( x \right){\zeta _n}\left( {p;y} \right) - {\rm{L}}{{\rm{i}}_p}\left( y \right){\zeta _n}\left( {m;x} \right)}}{n}{z^n}} \nonumber \\
  &= {\rm{L}}{{\rm{i}}_p}\left( y \right)\sum\limits_{n = 1}^\infty  {\frac{{{\zeta _n}\left( {1;z} \right)}}{{{n^m}}}{x^n}}  - {\rm{L}}{{\rm{i}}_m}\left( x \right)\sum\limits_{n = 1}^\infty  {\frac{{{\zeta _n}\left( {1;z} \right)}}{{{n^p}}}{y^n}}\nonumber\\&\quad+ {\rm{L}}{{\rm{i}}_m}\left( x \right){\rm{L}}{{\rm{i}}_{p + 1}}\left( {zy} \right) - {\rm{L}}{{\rm{i}}_p}\left( y \right){\rm{L}}{{\rm{i}}_{m + 1}}\left( {zx} \right),
\end{align}
where the partial sum ${{\zeta _n}\left( {m;a} \right)}$ is defined by \[{\zeta _n}\left( {m;a} \right) = \sum\limits_{k = 1}^n {\frac{{{a^k}}}{{{k^m}}}},\quad (a \in \mathbb{C}).\]
By convenience, we let ${{\zeta _0}\left( {m;a} \right)}=0$.
It is obvious that
\[{\zeta _n}\left( {m;1} \right) = H_n^{\left( m \right)}\quad{\rm and}\quad {\zeta _n}\left( {m; - 1} \right) =  - \bar H_n^{\left( m \right)}.\]
\end{thm}
\pf Let us consider the $N$-th partial sum
\begin{align*}
& {S_N}\left( {x,y,z} \right) = \sum\limits_{n = 1}^N {\frac{{{\rm{L}}{{\rm{i}}_m}\left( x \right){\zeta _n}\left( {p;y} \right) - {\rm{L}}{{\rm{i}}_p}\left( y \right){\zeta _n}\left( {m;x} \right)}}{n}{z^n}}  \\
  &= {\rm{L}}{{\rm{i}}_m}\left( x \right)\sum\limits_{n = 1}^N {\frac{{{\zeta _n}\left( {p;y} \right)}}{n}{z^n}}  - {\rm{L}}{{\rm{i}}_p}\left( y \right)\sum\limits_{n = 1}^N {\frac{{{\zeta _n}\left( {m;x} \right)}}{n}{z^n}}  \\
  &= {\rm{L}}{{\rm{i}}_m}\left( x \right)\sum\limits_{n = 1}^N {\sum\limits_{k = 1}^n {\frac{{{y^k}{z^n}}}{{n{k^p}}}} }  - {\rm{L}}{{\rm{i}}_p}\left( y \right)\sum\limits_{n = 1}^N {\sum\limits_{k = 1}^n {\frac{{{x^k}{z^n}}}{{n{k^m}}}} }  \\
  &= {\rm{L}}{{\rm{i}}_m}\left( x \right)\sum\limits_{k = 1}^N {\sum\limits_{n = k}^N {\frac{{{y^k}{z^n}}}{{n{k^p}}}} }  - {\rm{L}}{{\rm{i}}_p}\left( y \right)\sum\limits_{n = 1}^N {\sum\limits_{n = k}^N {\frac{{{x^k}{z^n}}}{{n{k^m}}}} }  \\
  &= {\rm{L}}{{\rm{i}}_m}\left( x \right)\sum\limits_{k = 1}^N {\frac{{{\zeta _N}\left( {1;z} \right) - {\zeta _{k - 1}}\left( {1;z} \right)}}{{{k^p}}}{y^k}}  - {\rm{L}}{{\rm{i}}_p}\left( y \right)\sum\limits_{k = 1}^N {\frac{{{\zeta _N}\left( {1;z} \right) - {\zeta _{k - 1}}\left( {1;z} \right)}}{{{k^m}}}{x^k}}  \\
  &= {\zeta _N}\left( {1;z} \right)\left( {{\rm{L}}{{\rm{i}}_m}\left( x \right){\zeta _N}\left( {p;y} \right) - {\rm{L}}{{\rm{i}}_p}\left( y \right){\zeta _N}\left( {m;x} \right)} \right) \\
 &\quad + {\rm{L}}{{\rm{i}}_p}\left( y \right)\sum\limits_{k = 1}^N {\frac{{{\zeta _{k - 1}}\left( {1;z} \right)}}{{{k^m}}}{x^k}}  - {\rm{L}}{{\rm{i}}_m}\left( x \right)\sum\limits_{k = 1}^N {\frac{{{\zeta _{k - 1}}\left( {1;z} \right)}}{{{k^p}}}{y^k}}.
\end{align*}
Then letting $N\to \infty $ in above equation, we get the required identity (\ref{c1}). The proof of Theorem \ref{thm3.2} is thus completed.\hfill$\square$\\
Taking $(x,y,z)\to (1,1,1),(1,-1,1)$ and $(-1,-1,1)$ in Theorem \ref{thm3.2}, we can get the following corollary.
\begin{cor} \label{cor3.3}
For positive integers $m$ and $p$, then
\begin{align}
 \sum\limits_{n = 1}^\infty  {\frac{{\zeta \left( m \right)H_n^{\left( p \right)} - \zeta \left( p \right)H_n^{\left( m \right)}}}{n}}  &= \zeta \left( p \right){S_{1,m}} - \zeta \left( m \right){S_{1,p}} \nonumber\\ &\quad+ \zeta \left( m \right)\zeta \left( {p + 1} \right) - \zeta \left( p \right)\zeta \left( {m + 1} \right),\:\:\left( {m,p \ge 2} \right) \label{c3}\\
 \sum\limits_{n = 1}^\infty  {\frac{{\zeta \left( m \right)\bar H_n^{\left( p \right)} - \bar \zeta \left( p \right)H_n^{\left( m \right)}}}{n}} & = \bar \zeta \left( p \right){S_{1,m}} - \zeta \left( m \right){S_{1,\bar p}} \nonumber\\ &\quad+ \zeta \left( m \right)\bar \zeta \left( {p + 1} \right) - \bar \zeta \left( p \right)\zeta \left( {m + 1} \right),\:\:\:\left( {m \ge 2} \right) \label{c4}\\
 \sum\limits_{n = 1}^\infty  {\frac{{\bar \zeta \left( m \right)\bar H_n^{\left( p \right)} - \bar \zeta \left( p \right)\bar H_n^{\left( m \right)}}}{n}}  &= \bar \zeta \left( p \right){S_{1,\bar m}} - \bar \zeta \left( m \right){S_{1,\bar p}} \nonumber\\ &\quad+ \bar \zeta \left( m \right)\bar \zeta \left( {p + 1} \right) - \bar \zeta \left( p \right)\bar \zeta \left( {m + 1} \right). \label{c5}
\end{align}
\end{cor}

The main Theorems of this section can be stated as follows.
\begin{thm}\label{thm3.4}
For positive integers $m,p$ and real $a,b,x$ with the conditions $\left| a \right| = \left| b \right| > 1$ and $\left| x \right| < {\left| a \right|^{ - 1}}$ or ${\left| b \right|^{ - 1}}$, we have
\begin{align}\label{c6}
&{\left( { - 1} \right)^p}\sum\limits_{n = 1}^\infty  {\frac{{{\zeta _n}\left( {1;ax} \right){\zeta _n}\left( {m;b} \right)}}{{{a^n}{n^p}}}}  - {\left( { - 1} \right)^m}\sum\limits_{n = 1}^\infty  {\frac{{{\zeta _n}\left( {1;bx} \right){\zeta _n}\left( {p;a} \right)}}{{{b^n}{n^m}}}}  \nonumber\\
 = &  \sum\limits_{i = 1}^{p - 1} {{{\left( { - 1} \right)}^{i - 1}}{\rm{L}}{{\rm{i}}_{p + 1 - i}}\left( {ax} \right)\sum\limits_{n = 1}^\infty  {\frac{{{\zeta _n}\left( {m;b} \right)}}{{{n^i}}}{x^n}} }  - \sum\limits_{i = 1}^{m - 1} {{{\left( { - 1} \right)}^{i - 1}}{\rm{L}}{{\rm{i}}_{m + 1 - i}}\left( {bx} \right)\sum\limits_{n = 1}^\infty  {\frac{{{\zeta _n}\left( {p;a} \right)}}{{{n^i}}}{x^n}} } \nonumber\\
 & + {\left( { - 1} \right)^p}\ln \left( {1 - ax} \right)\sum\limits_{n = 1}^\infty  {\frac{{{\zeta _n}\left( {m;b} \right)}}{{{n^p}}}\left( {{x^n} - {a^{ - n}}} \right)}  - {\left( { - 1} \right)^m}\ln \left( {1 - bx} \right)\sum\limits_{n = 1}^\infty  {\frac{{{\zeta _n}\left( {p;a} \right)}}{{{n^m}}}\left( {{x^n} - {b^{ - n}}} \right)}.
\end{align}
\end{thm}
\pf Applying the definition of polylogarithm function ${\rm Li}_p(x)$ and using the Cauchy product of power series, we can find that
\[\frac{{{\rm{L}}{{\rm{i}}_p}\left( {at} \right)}}{{1 - t}} = \sum\limits_{n = 1}^\infty  {\zeta_n(p;a){t^n}} ,\;t \in \left( { - 1,1} \right).\]
Then using above identity it is easily see that
\begin{align*}
\int\limits_0^x {\frac{{{\rm{L}}{{\rm{i}}_p}\left( {at} \right){\rm{L}}{{\rm{i}}_m}\left( {bt} \right)}}{{t\left( {1 - t} \right)}}} dt =\sum\limits_{n = 1}^\infty  {\zeta_n(p;a)\int\limits_0^x {{t^{n - 1}}{\rm{L}}{{\rm{i}}_m}\left( {bt} \right)} dt}
 = \sum\limits_{n = 1}^\infty  {\zeta_n(m;b)\int\limits_0^x {{t^{n - 1}}{\rm{L}}{{\rm{i}}_p}\left( {at} \right)} dt} .
\end{align*}
Then with the help of formula (\ref{c1}) and the elementary transform
\[\int\limits_0^x {{t^{n - 1}}{\rm{L}}{{\rm{i}}_m}\left( {bt} \right)dt}  = \frac{1}{{{b^n}}}\int\limits_0^{bx} {{t^{n - 1}}{\rm{L}}{{\rm{i}}_m}\left( t \right)dt} \]
by a direct calculation, we may deduce the result. \hfill$\square$

\begin{thm}\label{thm3.5}
For positive integers $m,p$ and real $a,b,x$ with the conditions $\left| a \right| = \left| b \right| > 1$ and $\left| x \right| < {\left| a \right|^{ - 1}}$ or ${\left| b \right|^{ - 1}}$, then the following identity holds
\begin{align}\label{c7}
&{\left( { - 1} \right)^p}\sum\limits_{n = 1}^\infty  {\frac{{{\zeta _n}\left( {2;ax} \right){\zeta _n}\left( {m;b} \right)}}{{{a^n}{n^p}}}}  - {\left( { - 1} \right)^m}\sum\limits_{n = 1}^\infty  {\frac{{{\zeta _n}\left( {2;bx} \right){\zeta _n}\left( {p;a} \right)}}{{{b^n}{n^m}}}}  \nonumber\\
 =& \sum\limits_{i = 1}^{p - 1} {\sum\limits_{j = 1}^{p - i} {{{\left( { - 1} \right)}^{i + j}}} {\rm{L}}{{\rm{i}}_{p + 2 - i - j}}\left( {ax} \right){\rm{L}}{{\rm{i}}_{i + j,m}}\left( {x,b} \right)} - \sum\limits_{i = 1}^{m - 1} {\sum\limits_{j = 1}^{m - i} {{{\left( { - 1} \right)}^{i + j}}} {\rm{L}}{{\rm{i}}_{m + 2 - i - j}}\left( {bx} \right){\rm{L}}{{\rm{i}}_{i + j,p}}\left( {x,a} \right)} \nonumber\\
 & + {\left( { - 1} \right)^p}p\ln \left( {1 - ax} \right)\left( {{\rm{L}}{{\rm{i}}_{p + 1,m}}\left( {x,b} \right) - {\rm{L}}{{\rm{i}}_{p + 1,m}}\left( {{a^{ - 1}},b} \right)} \right)\nonumber\\
 & - {\left( { - 1} \right)^m}m\ln \left( {1 - bx} \right)\left( {{\rm{L}}{{\rm{i}}_{m + 1,p}}\left( {x,a} \right) - {\rm{L}}{{\rm{i}}_{m + 1,p}}\left( {{b^{ - 1}},a} \right)} \right)\nonumber\\
 & - {\left( { - 1} \right)^p}p\sum\limits_{n = 1}^\infty  {\frac{{{\zeta _n}\left( {1;ax} \right){\zeta _n}\left( {m;b} \right)}}{{{a^n}{n^{p + 1}}}}}  + {\left( { - 1} \right)^m}m\sum\limits_{n = 1}^\infty  {\frac{{{\zeta _n}\left( {1;bx} \right){\zeta _n}\left( {p;a} \right)}}{{{b^n}{n^{m + 1}}}}} \nonumber \\
 & + {\left( { - 1} \right)^p}{\rm{L}}{{\rm{i}}_2}\left( {ax} \right){\rm{L}}{{\rm{i}}_{p,m}}\left( {{a^{ - 1}},b} \right) - {\left( { - 1} \right)^m}{\rm{L}}{{\rm{i}}_2}\left( {bx} \right){\rm{L}}{{\rm{i}}_{m,p}}\left( {{b^{ - 1}},a} \right),
\end{align}
where ${\rm{L}}{{\rm{i}}_{s,t}}\left( {x,y} \right)$ is defined by the double sum
\[{\rm{L}}{{\rm{i}}_{s,t}}\left( {x,y} \right): = \sum\limits_{n = 1}^\infty  {\frac{{{x^n}}}{{{n^s}}}\sum\limits_{k = 1}^n {\frac{{{y^k}}}{{{k^t}}}} } ,\;\left| {x} \right| \leq 1,\left| {y} \right| > 0,\left| {xy} \right| < 1,s,t>0\]
Of course, if $s>1$, then we can allow $\left| {xy} \right|=1$.
\end{thm}
\pf Replacing $x$ by $t$ in (\ref{c6}), then dividing it by $t$ and integrating over the interval $(0,x)$ with the help of formulae (\ref{c1}) and
\[\int\limits_0^x {\frac{{\ln \left( {1 - at} \right)}}{t}dt}  =  - {\rm{L}}{{\rm{i}}_2}\left( {ax} \right),\]
by a simple calculation, we have the result.\hfill$\square$\\
Letting $\left( {a,b,x} \right) \to \left( { - 1,1,1} \right),\left( { - 1, - 1,1} \right)$ and $\left( {p,a,b,x} \right) \to \left( {1, - 1,1, - 1} \right),\left( {1, - 1, - 1, - 1} \right)$ in Theorem \ref{thm3.4}, then with the help of Corollary \ref{cor3.3}  we obtain the following corollaries.
\begin{cor}\label{cor3.6}
For positive integers $p \geq 2 $ and $m \ge 2$, then
\begin{align}
{\left( { - 1} \right)^m}{S_{1\bar p,m}} + {\left( { - 1} \right)^p}{S_{m\bar 1,\bar p}} =& \sum\limits_{i = 2}^{m - 1} {{{\left( { - 1} \right)}^{i - 1}}\zeta \left( {m + 1 - i} \right){S_{\bar p,i}}}  - \sum\limits_{i = 2}^{p - 1} {{{\left( { - 1} \right)}^{i - 1}}\bar \zeta \left( {p + 1 - i} \right){S_{m,i}}}  \nonumber \\
&  + \bar \zeta \left( p \right){S_{1,m}} - \zeta \left( m \right){S_{1,\bar p}}\: + \zeta \left( m \right)\bar \zeta \left( {p + 1} \right) - \bar \zeta \left( p \right)\zeta \left( {m + 1} \right),\nonumber\\
&+ {\left( { - 1} \right)^p}\left( {{S_{m,p}} + {S_{m,\bar p}}} \right)\ln 2,\label{c8}\\
{\left( { - 1} \right)^p}{S_{\bar 1\bar m,\bar p}} - {\left( { - 1} \right)^m}{S_{\bar 1\bar p,\bar m}} =& \sum\limits_{i = 2}^{m - 1} {{{\left( { - 1} \right)}^{i - 1}}\bar \zeta \left( {m + 1 - i} \right){S_{\bar p,i}}}  - \sum\limits_{i = 2}^{p - 1} {{{\left( { - 1} \right)}^{i - 1}}\bar \zeta \left( {p + 1 - i} \right){S_{\bar m,i}}} \nonumber \\
& + \bar \zeta \left( p \right){S_{1,\bar m}} - \bar \zeta \left( m \right){S_{1,\bar p}} + \bar \zeta \left( m \right)\bar \zeta \left( {p + 1} \right) - \bar \zeta \left( p \right)\bar \zeta \left( {m + 1} \right)\nonumber \\
& + {\left( { - 1} \right)^p}\left( {{S_{\bar m,p}} + {S_{\bar m,\bar p}}} \right)\ln 2 - {\left( { - 1} \right)^m}\left( {{S_{\bar p,m}} + {S_{\bar p,\bar m}}} \right)\ln 2.\label{c9}
\end{align}
\end{cor}

\begin{cor}\label{cor3.7} For integer $m>1$, then
\begin{align}
&{\left( { - 1} \right)^{m - 1}}{S_{{{\bar 1}^2},m}} + {S_{1m,\bar 1}} = \sum\limits_{i = 1}^{m - 1} {{{\left( { - 1} \right)}^{i - 1}}\bar \zeta \left( {m + 1 - i} \right){S_{\bar 1,\bar i}}}  + {\left( { - 1} \right)^{m - 1}}\left( {{S_{\bar 1,m}} + {S_{\bar 1,\bar m}}} \right)\ln 2,\label{c10}\\
&{\left( { - 1} \right)^m}{S_{1\bar 1,\bar m}} + {S_{1\bar m,\bar 1}} = \sum\limits_{i = 1}^{m - 1} {{{\left( { - 1} \right)}^{i - 1}}\zeta \left( {m + 1 - i} \right){S_{\bar 1,\bar i}}} .\label{c11}
\end{align}
\end{cor}

Similarly, letting $(a,b,x)\rightarrow(1,1,1)$ and $ (-1,-1,-1)$ in Theorem \ref{thm3.5}, we can get the following two corollaries.

\begin{cor} \label{cor3.8}
 For positive integers $m>1$ and $p>1$, we have
\begin{align}\label{c12}
{\left( { - 1} \right)^p}{S_{2m,p}} - {\left( { - 1} \right)^m}{S_{2p,m}}\nonumber =& \sum\limits_{i = 1}^{p - 1} {\sum\limits_{j = 1}^{p - i} {{{\left( { - 1} \right)}^{i + j}}\zeta \left( {p + 2 - i - j} \right)} {S_{m,i + j}}}\nonumber \\
& - \sum\limits_{i = 1}^{m - 1} {\sum\limits_{j = 1}^{m - i} {{{\left( { - 1} \right)}^{i + j}}\zeta \left( {m + 2 - i - j} \right)} {S_{p,i + j}}} \nonumber\\
& - {\left( { - 1} \right)^p}p{S_{1m,p + 1}} + {\left( { - 1} \right)^m}m{S_{1p,m + 1}}\nonumber\\
& + {\left( { - 1} \right)^p}\zeta \left( 2 \right){S_{m,p}} - {\left( { - 1} \right)^m}\zeta \left( 2 \right){S_{p,m}}.
\end{align}
\end{cor}

\begin{cor}\label{cor3.9}
For positive integers $m$ and $p$, then
\begin{align}\label{c13}
 {\left( { - 1} \right)^m}{S_{2\bar p,\bar m}} - {\left( { - 1} \right)^p}{S_{2\bar m,\bar p}} =& \sum\limits_{i = 1}^{m - 1} {\sum\limits_{j = 1}^{m - i} {{{\left( { - 1} \right)}^{i + j}}\zeta \left( {m + 2 - i - j} \right){S_{\bar p,\overline {i + j}}}} } \nonumber  \\
  &- \sum\limits_{i = 1}^{p - 1} {\sum\limits_{j = 1}^{p - i} {{{\left( { - 1} \right)}^{i + j}}\zeta \left( {p + 2 - i - j} \right){S_{\bar m,\overline{i + j}}}} } \nonumber  \\
  &+ {\left( { - 1} \right)^{m - 1}}m{S_{1\bar p,\overline{m + 1}}} - {\left( { - 1} \right)^{p - 1}}p{S_{1\bar m,\overline{p + 1}}}\nonumber  \\
  &- {\left( { - 1} \right)^{m - 1}}\zeta \left( 2 \right){S_{\bar p,\bar m}} + {\left( { - 1} \right)^{p - 1}}\zeta \left( 2 \right){S_{\bar m,\bar p}}.
\end{align}
\end{cor}

\begin{thm}\label{thm3.10} For positive integer $k$ and $1 \le {i_1} \le {i_2} \le  \cdots  \le {i_k}$ with $x_j\in [-1,1]$, then the following relation holds:
\begin{align}\label{c14}
\prod\limits_{j = 1}^k {{\zeta _n}\left( {{i_j};{x_j}} \right)}  = \sum\limits_{j = 0}^{k - 1} {\sum\limits_{1 \le {l_1} <  \cdots  < {l_j} \le k} {{{\left( { - 1} \right)}^{k - j + 1}}\sum\limits_{m = 1}^n {\frac{{{\zeta _m}\left( {{i_{{l_1}}};{x_{{l_1}}}} \right) \cdots {\zeta _m}\left( {{i_{{l_j}}};{x_{{l_j}}}} \right)}}{{{m^{{i_1} +  \cdots  + {i_k} - \left( {{i_{{l_1}}} +  \cdots  + {i_{{l_j}}}} \right)}}}}\frac{{x_1^m \cdots x_k^m}}{{x_{{i_{{l_1}}}}^m \cdots x_{{i_{{l_j}}}}^m}}} } } ,
\end{align}
where by convention, the term corresponding to $j=0$ is ${\left( { - 1} \right)^{k + 1}}{\zeta _n}\left( {{i_1} +  \cdots  + {i_k};{x_1} \cdots {x_k}} \right)$.
\end{thm}
\pf  Consider the generating function
\[{F_k}\left( {x;{{\rm X}_k}} \right): = \sum\limits_{n = 1}^\infty  {\left\{ {\prod\limits_{j = 1}^k {{\zeta _n}\left( {{i_j};{x_j}} \right)}  - {\zeta _n}\left( {{i_1} +  \cdots  + {i_k};{x_1} \cdots {x_k}} \right)} \right\}{x^n}} ,\quad \left( {x \in \left( { - 1,1} \right)} \right),\]
where ${{\rm X}_k} = \left( {{x_1}, \cdots ,{x_k}} \right).$\\
Using the definition of the numbers ${{\zeta _n}\left( {m;a} \right)}$, we have
\begin{align}\label{c15}
 \left( {1 - x} \right){F_k}\left( {x;{{\rm X}_k}} \right) &= \sum\limits_{n = 1}^\infty  {\left\{ {\prod\limits_{j = 1}^k {{\zeta _n}\left( {{i_j};{x_j}} \right)}  - {\zeta _n}\left( {{i_1} +  \cdots  + {i_k};\prod\limits_{j = 1}^k {{x_j}} } \right)} \right\}{x^n}}  \nonumber\\
  &\quad- \sum\limits_{n = 1}^\infty  {\left\{ {\prod\limits_{j = 1}^k {{\zeta _{n - 1}}\left( {{i_j};{x_j}} \right)}  - {\zeta _{n - 1}}\left( {\sum\limits_{j = 1}^k {{i_j}} ;\prod\limits_{j = 1}^k {{x_j}} } \right)} \right\}{x^n}} \nonumber \\
  =& \sum\limits_{n = 1}^\infty  {\left\{ {\prod\limits_{j = 1}^k {{\zeta _n}\left( {{i_j};{x_j}} \right)}  - \prod\limits_{j = 1}^k {\left( {{\zeta _n}\left( {{i_j};{x_j}} \right) - \frac{{x_j^n}}{{{n^{{i_j}}}}}} \right)}  - \frac{{\prod\limits_{j = 1}^k {x_j^n} }}{{{n^{{i_1} +  \cdots  + {i_k}}}}}} \right\}{x^n}}\nonumber  \\
  =& - \sum\limits_{n = 1}^\infty  {\left\{ {\frac{{\prod\limits_{j = 1}^k {x_j^n} }}{{{n^{{i_1} +  \cdots  + {i_k}}}}} + \sum\limits_{j = 0}^{k - 1} {\sum\limits_{1 \le {l_1} <  \cdots  < {l_j} \le k} {{{\left( { - 1} \right)}^{k - j}}\frac{{\prod\limits_{r = 1}^j {{\zeta _n}\left( {{i_{{l_r}}};{x_{{l_r}}}} \right)} }}{{{n^{{i_1} +  \cdots  + {i_k} - \left( {{i_{{l_1}}} +  \cdots  + {i_{{l_j}}}} \right)}}}}\frac{{\prod\limits_{a = 1}^k {x_a^n} }}{{\prod\limits_{b = 1}^j {x_{{i_{{l_b}}}}^n} }}} } } \right\}{x^n}}
\end{align}
Multiplying this equation (\ref{c15}) by $1/(1-x)$, using the Cauchy product of power series, and identifying the coefficients of $x^n$, we may deduce the desired result. \hfill$\square$

\subsection{Identities for alternating Euler sums}
It is clear that from Theorem \ref{thm3.4}, \ref{thm3.5}, \ref{thm3.10} and Corollary \ref{cor3.6}-\ref{cor3.9}, we can establish many identities involving two or more Euler sums of the same weight. Some examples on alternating Euler sums with weight $\leq 5$ follows:
\begin{align}
&\sum\limits_{n = 1}^\infty  {\frac{{{H_n}{{\bar H}_n}}}{{{n^2}}}{{\left( { - 1} \right)}^{n - 1}}}  + \sum\limits_{n = 1}^\infty  {\frac{{{H_n}\bar H_n^{\left( 2 \right)}}}{n}{{\left( { - 1} \right)}^{n - 1}}}  = \frac{5}{4}\zeta \left( 4 \right) + \frac{1}{2}\zeta \left( 2 \right){\ln ^2}2,\label{c16}\\
&\sum\limits_{n = 1}^\infty  {\frac{{{H_n}H_n^{\left( 2 \right)}}}{n}{{\left( { - 1} \right)}^{n - 1}}}  - \sum\limits_{n = 1}^\infty  {\frac{{\bar H_n^2}}{{{n^2}}}}  = \frac{5}{8}\zeta \left( 4 \right) + \frac{7}{8}\zeta \left( 3 \right)\ln 2 - \frac{{11}}{4}\zeta \left( 2 \right){\ln ^2}2,\label{c17}\\
&\sum\limits_{n = 1}^\infty  {\frac{{{{\left( {H_n^{\left( 2 \right)}} \right)}^2}}}{n}{{\left( { - 1} \right)}^{n - 1}}}  + 2\sum\limits_{n = 1}^\infty  {\frac{{H_n^{\left( 2 \right)}{{\bar H}_n}}}{{{n^2}}}}  =  - \frac{{43}}{{16}}\zeta \left( 5 \right) + \frac{{15}}{8}\zeta \left( 2 \right)\zeta \left( 3 \right) + \frac{{35}}{8}\zeta \left( 4 \right)\ln 2, \label{c18}\\
& \sum\limits_{n = 1}^\infty  {\frac{{{{\left( {\bar H_n^{\left( 2 \right)}} \right)}^2}}}{n}{{\left( { - 1} \right)}^{n - 1}}}  + 2\sum\limits_{n = 1}^\infty  {\frac{{{{\bar H}_n}\bar H_n^{\left( 2 \right)}}}{{{n^2}}}{{\left( { - 1} \right)}^{n - 1}}}  = \frac{{19}}{{16}}\zeta \left( 5 \right) - \frac{1}{4}\zeta \left( 2 \right)\zeta \left( 3 \right) + \frac{5}{2}\zeta \left( 4 \right)\ln 2,\label{c19} \\
&\sum\limits_{n = 1}^\infty  {\frac{{{H_n}\bar H_n^{\left( 3 \right)}}}{n}{{\left( { - 1} \right)}^{n - 1}}}  - \sum\limits_{n = 1}^\infty  {\frac{{{H_n}{{\bar H}_n}}}{{{n^3}}}{{\left( { - 1} \right)}^{n - 1}}}  = \frac{1}{2}\zeta \left( 3 \right){\ln ^2}2 + \frac{9}{8}\zeta \left( 2 \right)\zeta \left( 3 \right) - \frac{{15}}{4}\zeta \left( 4 \right)\ln 2,\label{c20}\\
&\sum\limits_{n = 1}^\infty  {\frac{{{H_n}\bar H_n^{\left( 2 \right)}}}{{{n^2}}}}  + \sum\limits_{n = 1}^\infty  {\frac{{H_n^{\left( 2 \right)}{{\bar H}_n}}}{{{n^2}}}{{\left( { - 1} \right)}^{n - 1}}}  = 4{\rm{L}}{{\rm{i}}_4}\left( {\frac{1}{2}} \right)\ln 2 + \frac{1}{6}{\ln ^5}2 + \frac{7}{2}\zeta \left( 3 \right){\ln ^2}2 - \zeta \left( 2 \right){\ln ^3}2\nonumber\\&\quad\quad\quad\quad\quad\quad\quad\quad\quad\quad\quad\quad\quad\quad\quad\quad\quad + \frac{5}{8}\zeta \left( 2 \right)\zeta \left( 3 \right) - \frac{{23}}{{16}}\zeta \left( 4 \right)\ln 2,\label{c21}\\
&\sum\limits_{n = 1}^\infty  {\frac{{{{\bar H}_n}\bar H_n^{\left( 2 \right)}}}{{{n^2}}}}  + \sum\limits_{n = 1}^\infty  {\frac{{{H_n}H_n^{\left( 2 \right)}}}{{{n^2}}}{{\left( { - 1} \right)}^{n - 1}}}  =  - 4{\rm{L}}{{\rm{i}}_4}\left( {\frac{1}{2}} \right)\ln 2 - \frac{1}{6}{\ln ^5}2 - \frac{7}{2}\zeta \left( 3 \right){\ln ^2}2 + \zeta \left( 2 \right){\ln ^3}2 \nonumber\\&\quad\quad\quad\quad\quad\quad\quad\quad\quad\quad\quad\quad\quad\quad\quad\quad\quad+ \frac{3}{{16}}\zeta \left( 2 \right)\zeta \left( 3 \right) + \frac{{49}}{8}\zeta \left( 4 \right)\ln 2,\label{c22}\\
&\sum\limits_{n = 1}^\infty  {\frac{{{H_n}H_n^{\left( 3 \right)}}}{n}{{\left( { - 1} \right)}^{n - 1}}}  + \sum\limits_{n = 1}^\infty  {\frac{{\bar H_n^2}}{{{n^3}}}}  =  - 2{\rm{L}}{{\rm{i}}_4}\left( {\frac{1}{2}} \right)\ln 2 - \frac{1}{{12}}{\ln ^5}2 + \frac{{17}}{8}\zeta \left( 3 \right){\ln ^2}2 + \frac{1}{2}\zeta \left( 2 \right){\ln ^3}2 \nonumber\\&\quad\quad\quad\quad\quad\quad\quad\quad\quad\quad\quad\quad\quad\quad\quad+ \frac{{11}}{{16}}\zeta \left( 2 \right)\zeta \left( 3 \right) - \frac{{11}}{{16}}\zeta \left( 4 \right)\ln 2,\label{c23}\\
&\sum\limits_{n = 1}^\infty  {\frac{{H_n^{\left( 2 \right)}{{\bar H}_n}}}{{{n^2}}}{{\left( { - 1} \right)}^{n - 1}}}  + \sum\limits_{n = 1}^\infty  {\frac{{H_n^{\left( 2 \right)}\bar H_n^{\left( 2 \right)}}}{n}{{\left( { - 1} \right)}^{n - 1}}}  + \sum\limits_{n = 1}^\infty  {\frac{{{H_n}\bar H_n^{\left( 2 \right)}}}{{{n^2}}}{{\left( { - 1} \right)}^{n - 1}}}\nonumber\\& =  - 4{\rm{L}}{{\rm{i}}_5}\left( {\frac{1}{2}} \right) + \frac{1}{{30}}{\ln ^5}2 + \frac{{37}}{8}\zeta \left( 5 \right) + \frac{7}{4}\zeta \left( 3 \right){\ln ^2}2 - \frac{1}{3}\zeta \left( 2 \right){\ln ^3}2 + \frac{5}{8}\zeta \left( 2 \right)\zeta \left( 3 \right) - 3\zeta \left( 4 \right)\ln 2.\label{c24}
\end{align}
These identities can be obtained from the Theorem \ref{thm3.4}, \ref{thm3.5}, \ref{thm3.10} and Corollary \ref{cor3.6}-\ref{cor3.9} which are presented in the section \ref{sec3.1}. For example, taking $k=3,\ (i_1,i_2,i_3)=(2,2,1)$ and $(x_1,x_2,x_3)=(1,1,-1)$ or $(-1,-1,-1)$ in Theorem \ref{thm3.10}, then letting $n\rightarrow \infty$ yields the formulas (\ref{c18}) and (\ref{c19}).

In (\ref{b5}), setting $k=3,4$ and $x=-1$, using (\ref{b18}) we can get the following two identities

\begin{align}
&\sum\limits_{n = 1}^\infty  {\frac{{H_n^3 - 3{H_n}H_n^{\left( 2 \right)} + 2H_n^{\left( 3 \right)}}}{n}{{\left( { - 1} \right)}^{n - 1}}}  = 6\zeta \left( 4 \right) + \frac{3}{2}\zeta \left( 2 \right){\ln ^2}2 - \frac{1}{2}{\ln ^4}2 - \frac{{21}}{4}\zeta \left( 3 \right)\ln 2 - 6{\rm{L}}{{\rm{i}}_4}\left( {\frac{1}{2}} \right),\label{c25}\\
&\sum\limits_{n = 1}^\infty  {\frac{{H_n^4 + 8{H_n}H_n^{\left( 3 \right)} - 6H_n^2H_n^{\left( 2 \right)} + 3{{\left( {H_n^{\left( 2 \right)}} \right)}^2} - 6H_n^{\left( 4 \right)}}}{n}{{\left( { - 1} \right)}^{n - 1}}} \nonumber\\ &= 24{\rm{L}}{{\rm{i}}_5}\left( {\frac{1}{2}} \right) + 24{\rm{L}}{{\rm{i}}_4}\left( {\frac{1}{2}} \right)\ln 2 + {\ln ^5}2 + \frac{{21}}{2}\zeta \left( 3 \right){\ln ^2}2 - 24\zeta \left( 5 \right) - 4\zeta \left( 2 \right){\ln ^3}2.\label{c26}
\end{align}

\subsection{Closed form of alternating Euler sums}
In this subsection, we use our equation system to obtain explicit evaluations for alternating Euler sums with weight $\leq 5$.

Note that the linear sums $S_{\bar p,q},\ S_{p,\bar q}$ and $S_{\bar p,\bar q}$ are reducible to zeta values and polylogarithms when
$p + q$ is odd (see \cite{BBG1995,C2010,FS1998,Xu2016}). The quadratic sums ${S_{{1^2},{\overline {2m}}}},\ {S_{{{\bar 1}^2},{\overline {2m}}}} ,\ {S_{1\bar 1,\overline {2m}}} ,
{S_{1\bar 1,2m}} ,\ {S_{{{\bar 1}^2},2m}} $ can be computed by the paper \cite{Xyz2016}. The quadratic and cubic sums $S_{1^2,\bar 3},\ S_{{\bar 1}^2,3},\ S_{{\bar 1}^2,\bar 3},\ S_{1^3,\bar 2},\ S_{{\bar 1}^3,2}$ and $S_{{\bar 1}^3,\bar 2}$ can be evaluated by Bailey, Borwein and Girgensohn \cite{BBG1994}, we get
\begin{align}
& \sum\limits_{n = 1}^\infty  {\frac{{H_n^3}}{{{n^2}}}} {\left( { - 1} \right)^{n - 1}} = 6{\rm{L}}{{\rm{i}}_5}\left( {\frac{1}{2}} \right) + 6{\rm{L}}{{\rm{i}}_4}\left( {\frac{1}{2}} \right)\ln 2 + \frac{1}{5}{\ln ^5}2 + \frac{{21}}{8}\zeta \left( 3 \right){\ln ^2}2\nonumber \\
 &\quad\quad\quad\quad\quad\quad\quad\quad- \frac{9}{4}\zeta \left( 5 \right) - \zeta \left( 2 \right){\ln ^3}2 - \frac{{27}}{{16}}\zeta \left( 2 \right)\zeta \left( 3 \right), \label{c27}\\
& \sum\limits_{n = 1}^\infty  {\frac{{H_n^2}}{{{n^3}}}} {\left( { - 1} \right)^{n - 1}} = 4{\rm{L}}{{\rm{i}}_5}\left( {\frac{1}{2}} \right) + 4{\rm{L}}{{\rm{i}}_4}\left( {\frac{1}{2}} \right)\ln 2 + \frac{2}{{15}}{\ln ^5}2 + \frac{7}{4}\zeta \left( 3 \right){\ln ^2}2\nonumber\\
  &\quad\quad\quad\quad\quad\quad\quad\quad- \frac{{19}}{{32}}\zeta \left( 5 \right) - \frac{2}{3}\zeta \left( 2 \right){\ln ^3}2 - \frac{{11}}{8}\zeta \left( 2 \right)\zeta \left( 3 \right),\label{c28}\\
&\sum\limits_{n = 1}^\infty  {\frac{{\bar H_n^2}}{{{n^3}}}}= 4{\rm{L}}{{\rm{i}}_5}\left( {\frac{1}{2}} \right) - \frac{1}{{30}}{\ln ^5}2 + \frac{7}{4}\zeta \left( 3 \right){\ln ^2}2 - \frac{{167}}{{32}}\zeta \left( 5 \right) + \frac{1}{3}\zeta \left( 2 \right){\ln ^3}2\nonumber
\\&\quad\quad\quad\quad\quad\quad+ \frac{3}{4}\zeta \left( 2 \right)\zeta \left( 3 \right) + \frac{{19}}{8}\zeta \left( 4 \right)\ln 2,\label{c29}\\
&\sum\limits_{n = 1}^\infty  {\frac{{\bar H_n^2}}{{{n^3}}}{{\left( { - 1} \right)}^{n - 1}}} =  - 4{\rm{L}}{{\rm{i}}_4}\left( {\frac{1}{2}} \right)\ln 2 - \frac{1}{6}{\ln ^5}2 - \frac{{19}}{{32}}\zeta \left( 5 \right) + \zeta \left( 2 \right){\ln ^3}2 \nonumber\\&\quad\quad\quad\quad\quad\quad\quad\quad\quad\quad+ \frac{3}{8}\zeta \left( 2 \right)\zeta \left( 3 \right) + \frac{{19}}{8}\zeta \left( 4 \right)\ln 2.\label{c30}
 \end{align}
Hence, from identities (\ref{b46})-(\ref{b49}), (\ref{c16}), (\ref{c17}) and (\ref{c25}), we deduce the following six closed form
\begin{align}
&\sum\limits_{n = 1}^\infty  {\frac{{H_n^{\left( 2 \right)}{{\bar H}_n}}}{n}{{\left( { - 1} \right)}^{n - 1}}} = \frac{{43}}{{16}}\zeta \left( 4 \right) - \frac{1}{4}\zeta \left( 2 \right){\ln ^2}2 - \frac{1}{8}{\ln ^4}2 - 3{\rm{L}}{{\rm{i}}_4}\left( {\frac{1}{2}} \right),\label{c31}\\
&\sum\limits_{n = 1}^\infty  {\frac{{H_n^2{{\bar H}_n}}}{n}} {\left( { - 1} \right)^{n - 1}}= \frac{{53}}{{16}}\zeta \left( 4 \right) + \frac{1}{8}{\ln ^4}2 + \frac{1}{4}\zeta \left( 2 \right){\ln ^2}2 - 3{\rm{L}}{{\rm{i}}_4}\left( {\frac{1}{2}} \right),\label{c32}\\
&\sum\limits_{n = 1}^\infty  {\frac{{\bar H_n^3}}{n}{{\left( { - 1} \right)}^{n - 1}}} = 5{\rm{L}}{{\rm{i}}_4}\left( {\frac{1}{2}} \right) - \frac{{59}}{{16}}\zeta \left( 4 \right) + \frac{{13}}{4}\zeta \left( 2 \right){\ln ^2}2 + \frac{{11}}{{24}}{\ln ^4}2,\label{c33}\\
&\sum\limits_{n = 1}^\infty  {\frac{{H_n^{}\bar H_n^{\left( 2 \right)}}}{n}} {\left( { - 1} \right)^{n - 1}}= 3{\rm{L}}{{\rm{i}}_4}\left( {\frac{1}{2}} \right) + \frac{1}{8}{\ln ^4}2 - \frac{1}{4}\zeta \left( 2 \right){\ln ^2}2 - \frac{9}{{16}}\zeta \left( 4 \right),\label{c34}\\
&\sum\limits_{n = 1}^\infty  {\frac{{H_n^{}H_n^{\left( 2 \right)}}}{n}} {\left( { - 1} \right)^{n - 1}}= 2{\rm{L}}{{\rm{i}}_4}\left( {\frac{1}{2}} \right) + \frac{1}{{12}}{\ln ^4}2{\rm{ + }}\frac{7}{8}\zeta \left( 3 \right)\ln 2 - \frac{1}{4}\zeta \left( 2 \right){\ln ^2}2 - \zeta \left( 4 \right),\label{c35}\\
&\sum\limits_{n = 1}^\infty  {\frac{{H_n^3}}{n}{{\left( { - 1} \right)}^{n - 1}}}  = \frac{5}{8}\zeta \left( 4 \right) + \frac{3}{4}\zeta \left( 2 \right){\ln ^2}2 - \frac{1}{4}{\ln ^4}2 - \frac{9}{8}\zeta \left( 3 \right)\ln 2.\label{c36}
 \end{align}
Note that formula (\ref{c36}) was also proved in \cite{M2014} by another method. \\
In (2.5) and (3.8) of \cite{Xyz2016}, letting $s=3$ and $(p,l)=(1,3)$, respectively, and combining (\ref{c28})-(\ref{c30}) of this paper, we give
\begin{align}
\sum\limits_{n = 1}^\infty  {\frac{{{H_n}{{\bar H}_n}}}{{{n^3}}}{{\left( { - 1} \right)}^{n - 1}}}  =& 2{\rm{L}}{{\rm{i}}_5}\left( {\frac{1}{2}} \right) - \frac{1}{{60}}{\ln ^5}2 - \frac{{37}}{{16}}\zeta \left( 5 \right) - \frac{7}{8}\zeta \left( 3 \right){\ln ^2}2\nonumber  \\&+ \frac{1}{6}\zeta \left( 2 \right){\ln ^3}2+ 4\zeta \left( 4 \right)\ln 2 - \frac{1}{8}\zeta \left( 2 \right)\zeta \left( 3 \right),\label{c37}\\
\sum\limits_{n = 1}^\infty  {\frac{{{H_n}{{\bar H}_n}}}{{{n^3}}}} =&2{\rm Li}{_5}\left( {\frac{1}{2}} \right) - \frac{1}{{60}}{\ln ^5}2 - \frac{{193}}{{64}}\zeta \left( 5 \right) - \frac{7}{8}\zeta \left( 3 \right){\ln ^2}2\nonumber  \\
 &+ \frac{1}{6}\zeta \left( 2 \right){\ln ^3}2 + 4\zeta \left( 4 \right)\ln 2  + \frac{3}{8}\zeta \left( 2 \right)\zeta \left( 3 \right).\label{c38}
\end{align}
Combining (\ref{c20}) with (\ref{c37}), and (\ref{c23}) with (\ref{c29}), we have the two results
\begin{align}
\sum\limits_{n = 1}^\infty  {\frac{{{H_n}\bar H_n^{\left( 3 \right)}}}{n}} {\left( { - 1} \right)^{n - 1}}=& 2{\rm{L}}{{\rm{i}}_5}\left( {\frac{1}{2}} \right) - \frac{1}{{60}}{\ln ^5}2 - \frac{3}{8}\zeta \left( 3 \right){\ln ^2}2 - \frac{{37}}{{16}}\zeta \left( 5 \right) + \frac{1}{6}\zeta \left( 2 \right){\ln ^3}2\nonumber
 \\&+ \zeta \left( 2 \right)\zeta \left( 3 \right) + \frac{1}{4}\zeta \left( 4 \right)\ln 2,\label{c39}\\
\sum\limits_{n = 1}^\infty  {\frac{{{H_n}H_n^{\left( 3 \right)}}}{n}{{\left( { - 1} \right)}^{n - 1}}}  = & - 4{\rm{L}}{{\rm{i}}_5}\left( {\frac{1}{2}} \right) - 2{\rm{L}}{{\rm{i}}_4}\left( {\frac{1}{2}} \right)\ln 2 - \frac{1}{{20}}{\ln ^5}2 + \frac{3}{8}\zeta \left( 3 \right){\ln ^2}2 + \frac{{167}}{{32}}\zeta \left( 5 \right)
\nonumber  \\&+ \frac{1}{6}\zeta \left( 2 \right){\ln ^3}2 - \frac{1}{{16}}\zeta \left( 2 \right)\zeta \left( 3 \right) - \frac{{49}}{{16}}\zeta \left( 4 \right)\ln 2.\label{c40}
\end{align}
By (\ref{b18}), letting $m=4,z=1$, we arrive at the conclusion that
\begin{align}\label{c41}
\int\limits_0^1 {\frac{{{{\ln }^4}\left( {1 + x} \right)}}{x}} dx
 &=4\sum\limits_{n = 1}^\infty  {\frac{{H_n^3 - 3{H_n}H_n^{\left( 2 \right)} + 2H_n^{\left( 3 \right)}}}{{{{\left( {n + 1} \right)}^2}}}{{\left( { - 1} \right)}^{n + 1}}} \nonumber\\
           &=- 24{\rm Li}{_5}\left( {\frac{1}{2}} \right) - 24{\rm Li}{_4}\left( {\frac{1}{2}} \right)\ln 2 - \frac{4}{5}{\ln ^5}2 - \frac{{21}}{2}\zeta \left( 3 \right){\ln ^2}2 + 24\zeta \left( 5 \right) + 4\zeta \left( 2 \right){\ln ^3}2.
\end{align}
Substituting (\ref{c27}) into (\ref{c41}) yields
\begin{align}\label{c42}
\sum\limits_{n = 1}^\infty  {\frac{{{H_n}H_n^{\left( 2 \right)}}}{{{n^2}}}{{\left( { - 1} \right)}^{n - 1}}}
 &= - 4{\rm Li}{_5}\left( {\frac{1}{2}} \right) - 4{\rm Li}{_4}\left( {\frac{1}{2}} \right)\ln 2 - \frac{2}{{15}}{\ln ^5}2 - \frac{7}{4}\zeta \left( 3 \right){\ln ^2}2\nonumber\\
           &\quad \ + \frac{{23}}{8}\zeta \left( 5 \right) + \frac{2}{3}\zeta \left( 2 \right){\ln ^3}2 + \frac{{15}}{{16}}\zeta \left( 2 \right)\zeta \left( 3 \right).
\end{align}
Combining (\ref{c22}) with (\ref{c42}) gives
\begin{align}\label{c43}
\sum\limits_{n = 1}^\infty  {\frac{{{{\bar H}_n}\bar H_n^{\left( 2 \right)}}}{{{n^2}}}}
 &= 4{\rm Li}{_5}\left( {\frac{1}{2}} \right) - \frac{1}{{30}}{\ln ^5}2 - \frac{7}{4}\zeta \left( 3 \right){\ln ^2}2 - \frac{{23}}{8}\zeta \left( 5 \right) + \frac{1}{3}\zeta \left( 2 \right){\ln ^3}2
\nonumber \\
           &\quad \ - \frac{3}{4}\zeta \left( 2 \right)\zeta \left( 3 \right) + \frac{{49}}{8}\zeta \left( 4 \right)\ln 2.
\end{align}

\begin{thm}
For positive integer $m$, then the integrals
\[\int\limits_0^1 {\frac{{{{\ln }^{2m}}x{{\ln }^2}\left( {1 - x} \right)}}{{1 + x}}} dx\quad {\rm and}\quad \int\limits_0^1 {\frac{{{{\ln }^{2m}}x{{\ln }^2}\left( {1 + x} \right)}}{{1 - x}}} dx\]
are expressible in terms of $\ln2$, zeta values and linear sums.
\end{thm}
\pf Using (\ref{b32}) and the well known identity
\[\int\limits_0^1 {\frac{{{x^n}{{\ln }^m}x}}{{1 + x}}} dx = {\left( { - 1} \right)^{m + n}}m!\left( {\bar \zeta \left( {m + 1} \right) - \bar H_n^{\left( {m + 1} \right)}} \right),\quad m \in \N,\]
we can find that
\begin{align}
\int\limits_0^1 {\frac{{{{\ln }^m}x{{\ln }^2}\left( {1 - x} \right)}}{{1 + x}}} dx &= 2{\left( { - 1} \right)^m}m!\left\{ {\sum\limits_{n = 1}^\infty  {\frac{{{H_n}\bar H_n^{\left( {m + 1} \right)}}}{n}} {{\left( { - 1} \right)}^{n - 1}} - \sum\limits_{n = 1}^\infty  {\frac{{\bar H_n^{\left( {m + 1} \right)}}}{{{n^2}}}} {{\left( { - 1} \right)}^{n - 1}}} \right\}\nonumber\\
&\quad+{\left( { - 1} \right)^m}m!\bar \zeta \left( {m + 1} \right){\ln ^2}2,\label{c44}\\
\int\limits_0^1 {\frac{{{{\ln }^m}x{{\ln }^2}\left( {1 + x} \right)}}{{1 - x}}} dx &= 2{\left( { - 1} \right)^m}m!\left\{ {\sum\limits_{n = 1}^\infty  {\frac{{{H_n}H_n^{\left( {m + 1} \right)}}}{n}} {{\left( { - 1} \right)}^{n - 1}} - \sum\limits_{n = 1}^\infty  {\frac{{H_n^{\left( {m + 1} \right)}}}{{{n^2}}}} {{\left( { - 1} \right)}^{n - 1}}} \right\}\nonumber\\&\quad +{\left( { - 1} \right)^m}m!\zeta \left( {m + 1} \right){\ln ^2}2.\label{c45}
\end{align}
From Corollary \ref{cor3.7}, we see that the quadratic sums $S_{1(2m),\bar 1}$ and $S_{1 \overline {2m},\bar 1}$ are reducible to linear sums and zeta values. Thus, Theorem holds.\hfill$\square$\\
Setting $m=2$ in (\ref{c44}) and (\ref{c45}), and combining (\ref{c39}) and (\ref{c40}), we obtain
\begin{align}
\int\limits_0^1 {\frac{{{{\ln }^2}x{{\ln }^2}\left( {1 - x} \right)}}{{1 + x}}} dx &= 8{\rm{L}}{{\rm{i}}_5}\left( {\frac{1}{2}} \right) - \frac{1}{{15}}{\ln ^5}2 + \frac{{15}}{2}\zeta \left( 5 \right) + \frac{2}{3}\zeta \left( 2 \right){\ln ^3}2 - \frac{{13}}{2}\zeta \left( 2 \right)\zeta \left( 3 \right) + \zeta \left( 4 \right)\ln 2,\label{c46}\\
\int\limits_0^1 {\frac{{{{\ln }^2}x{{\ln }^2}\left( {1 + x} \right)}}{{1 - x}}} dx
 &= - 16{\rm Li}{_5}\left( {\frac{1}{2}} \right) - 8\ln 2{\rm Li}{_4}\left( {\frac{1}{2}} \right) - \frac{1}{5}{\ln ^5}2 + \frac{7}{2}\zeta \left( 3 \right){\ln ^2}2 + \frac{{47}}{2}\zeta \left( 5 \right)\nonumber \\
           &\quad \ + \frac{2}{3}\zeta \left( 2 \right){\ln ^3}2 - \frac{{13}}{4}\zeta \left( 2 \right)\zeta \left( 3 \right) - \frac{{49}}{4}\zeta \left( 4 \right)\ln 2.\label{c47}
\end{align}
Now, we consider the generating function \[y = \sum\limits_{n = 1}^\infty  {\bar H_n^2} {x^n},\quad( - 1 < x < 1).\]
By a direct calculation, it is easy to see that
\begin{align}\label{c48}
y = xy + 2\int\limits_0^x {\frac{{\ln \left( {1 - t} \right)}}{{1 + t}}dt}  + \sum\limits_{n = 1}^\infty  {\frac{{{x^n}}}{{{n^2}}}} .
\end{align}
By integrating, we may rewrite (\ref{c48}) as
\begin{align}\label{c49}
\sum\limits_{n = 1}^\infty  {\frac{{{H^{(2)} _n} - \bar H_n^2}}{{n + 1}}} {x^{n + 1}} = 2\ln \left( {1 - x} \right)\int_0^x {\frac{{\ln \left( {1 - t} \right)}}{{1 + t}}} dt - 2\int_0^x {\frac{{{{\ln }^2}\left( {1 - t} \right)}}{{1 + t}}} dt.
\end{align}
Multiplying (\ref{c49}) by $\frac{\ln{x}}{x}$ and integrating over the interval (0,1). The result is
\begin{align}\label{c50}
\sum\limits_{n = 1}^\infty  {\frac{{H_n^{\left( 2 \right)}{{\bar H}_n}}}{{{n^2}}}} {\left( { - 1} \right)^{n - 1}}
 &= \frac{5}{4}\zeta \left( 5 \right) + \frac{{15}}{4}\zeta \left( 4 \right)\ln 2 - \frac{9}{8}\zeta \left( 2 \right)\zeta \left( 3 \right) + \frac{1}{2}\sum\limits_{n = 1}^\infty  {\frac{{\bar H_n^2}}{{{{\left( {n + 1} \right)}^3}}}}
\nonumber \\
           &\quad \ - \sum\limits_{n = 1}^\infty  {\frac{{{H_n}{{\bar H}_n}}}{{{n^3}}}{{\left( { - 1} \right)}^{n - 1}} - \frac{1}{2}} \int\limits_0^1 {\frac{{{{\ln }^2}x{{\ln }^2}\left( {1 - x} \right)}}{{1 + x}}} dx.
\end{align}
Substituting (\ref{c29}), (\ref{c37}) and (\ref{c46}) into (\ref{c50}), we have
\begin{align}\label{c51}
\sum\limits_{n = 1}^\infty  {\frac{{H_n^{\left( 2 \right)}{{\bar H}_n}}}{{{n^2}}}} {\left( { - 1} \right)^{n - 1}}
 &= - 4{\rm Li}{_5}\left( {\frac{1}{2}} \right) + \frac{1}{{30}}{\ln ^5}2 + \frac{7}{4}\zeta \left( 3 \right){\ln ^2}2 - \frac{{29}}{{64}}\zeta \left( 5 \right)
\nonumber \\
           &\quad \ - \frac{1}{3}\zeta \left( 2 \right){\ln ^3}2+ \frac{{15}}{8}\zeta \left( 2 \right)\zeta \left( 3 \right) - \frac{{23}}{{16}}\zeta \left( 4 \right)\ln 2.
\end{align}
Applying the same argument as above, replacing $x$ by $-x$ in (\ref{c49}), then dividing this by $\frac{\ln{x}}{x}$ and integrating over the interval $(0,1)$, we arrive at the conclusion that
\begin{align}\label{c52}
\sum\limits_{n = 1}^\infty  {\frac{{{{\bar H}_n}\bar H_n^{\left( 2 \right)}}}{{{n^2}}}} {\left( { - 1} \right)^{n - 1}}
 &= 8{\rm Li}{_5}\left( {\frac{1}{2}} \right) + 4{\rm Li}{_4}\left( {\frac{1}{2}} \right)\ln 2 + \frac{1}{{10}}{\ln ^5}2 - \frac{{35}}{4}\zeta \left( 5 \right)\nonumber\\
           &\quad \ - \frac{1}{3}\zeta \left( 2 \right){\ln ^3}2+ \frac{1}{{16}}\zeta \left( 2 \right)\zeta \left( 3 \right) + \frac{{49}}{8}\zeta \left( 4 \right)\ln 2.
\end{align}
Letting $(a,x,b,p,m)=(1,1,-1,2,1),\ (-1,-1,1,1,2),\ (-1,1,-1,1,2)$ and $(-1,-1,-1,2,1)$ in Theorem \ref{thm3.5}, and combining identities (\ref{b51}), (\ref{c18}), (\ref{c19}), (\ref{c21}), (\ref{c26}), (\ref{c37}), (\ref{c38}), (\ref{c40}), (\ref{c42}), (\ref{c43}), (\ref{c51}) and (\ref{c52}), we obtain the following eight evaluations
\begin{align}
&\sum\limits_{n = 1}^\infty  {\frac{{{H_n}\bar H_n^{\left( 2 \right)}}}{{{n^2}}}} {\left( { - 1} \right)^{n - 1}} = \frac{{125}}{{32}}\zeta \left( 5 \right) - \frac{{13}}{8}\zeta \left( 2 \right)\zeta \left( 3 \right),\label{c53}\\
&\sum\limits_{n = 1}^\infty  {\frac{{H_n^{\left( 2 \right)}\bar H_n^{\left( 2 \right)}}}{n}{{\left( { - 1} \right)}^{n - 1}}}  = \frac{{75}}{{64}}\zeta \left( 5 \right) + \frac{3}{8}\zeta \left( 2 \right)\zeta \left( 3 \right) - \frac{{25}}{{16}}\zeta \left( 4 \right)\ln 2,\label{c54}\\
&\sum\limits_{n = 1}^\infty  {\frac{{{H_n}\bar H_n^{\left( 2 \right)}}}{{{n^2}}}}  = 4{\rm{L}}{{\rm{i}}_5}\left( {\frac{1}{2}} \right) + 4{\rm{L}}{{\rm{i}}_4}\left( {\frac{1}{2}} \right)\ln 2 + \frac{2}{{15}}{\ln ^5}2 + \frac{{29}}{{64}}\zeta \left( 5 \right) + \frac{7}{4}\zeta \left( 3 \right){\ln ^2}2\nonumber\\&\quad \quad \quad \quad \quad \quad \quad - \frac{2}{3}\zeta \left( 2 \right){\ln ^3}2 - \frac{5}{4}\zeta \left( 2 \right)\zeta \left( 3 \right),\label{c55}\\
&\sum\limits_{n = 1}^\infty  {\frac{{H_n^2H_n^{\left( 2 \right)}}}{n}} {\left( { - 1} \right)^{n - 1}} =  - 2{\rm{L}}{{\rm{i}}_4}\left( {\frac{1}{2}} \right)\ln 2 - \frac{1}{{12}}{\ln ^5}2 - \frac{7}{8}\zeta \left( 3 \right){\ln ^2}2 + \frac{1}{3}\zeta \left( 2 \right){\ln ^3}2\nonumber\\
&\quad\quad\quad\quad\quad\quad\quad\quad\quad\quad + \frac{3}{8}\zeta \left( 2 \right)\zeta \left( 3 \right) + \frac{7}{8}\zeta \left( 4 \right)\ln 2,\label{c56}\\
&\sum\limits_{n = 1}^\infty  {\frac{{H_n^{(2)}{{\bar H}_n}}}{{{n^2}}}}  =  - 8{\rm{L}}{{\rm{i}}_5}\left( {\frac{1}{2}} \right) - 4{\rm{L}}{{\rm{i}}_4}\left( {\frac{1}{2}} \right)\ln 2 - \frac{1}{{10}}{\ln ^5}2 + \frac{{27}}{4}\zeta \left( 5 \right) + \frac{1}{3}\zeta \left( 2 \right){\ln ^3}2\nonumber\\
&\quad\quad\quad\quad\quad\quad + \frac{1}{2}\zeta \left( 2 \right)\zeta \left( 3 \right) - \frac{{23}}{{16}}\zeta \left( 4 \right)\ln 2,\label{c57}\\
&\sum\limits_{n = 1}^\infty  {\frac{{{{\left( {H_n^{\left( 2 \right)}} \right)}^2}}}{n}} {\left( { - 1} \right)^{n - 1}} = 16{\rm{L}}{{\rm{i}}_5}\left( {\frac{1}{2}} \right) + 8{\rm{L}}{{\rm{i}}_4}\left( {\frac{1}{2}} \right)\ln 2 + \frac{1}{5}{\ln ^5}2 - \frac{{259}}{{16}}\zeta \left( 5 \right) - \frac{2}{3}\zeta \left( 2 \right){\ln ^3}2\nonumber\\
&\quad\quad\quad\quad\quad\quad\quad\quad\quad\quad +\frac{5}{8}\zeta \left( 2 \right)\zeta \left( 3 \right) + \frac{{29}}{4}\zeta \left( 4 \right)\ln 2,\label{c58}\\
&\sum\limits_{n = 1}^\infty  {\frac{{{{\left( {\bar H_n^{\left( 2 \right)}} \right)}^2}}}{n}{{\left( { - 1} \right)}^{n - 1}}}  =  - 16{\rm{L}}{{\rm{i}}_5}\left( {\frac{1}{2}} \right) - 8{\rm{L}}{{\rm{i}}_4}\left( {\frac{1}{2}} \right)\ln 2 - \frac{1}{5}{\ln ^5}2 + \frac{{299}}{{16}}\zeta \left( 5 \right) + \frac{2}{3}\zeta \left( 2 \right){\ln ^3}2\nonumber\\
&\quad\quad\quad\quad\quad\quad\quad\quad\quad\quad - \frac{1}{4}\zeta \left( 2 \right)\zeta \left( 3 \right) - \frac{{39}}{4}\zeta \left( 4 \right)\ln 2,\label{c59}\\
&\sum\limits_{n = 1}^\infty  {\frac{{H_n^4}}{n}} {\left( { - 1} \right)^{n - 1}} = 8{\rm{L}}{{\rm{i}}_5}\left( {\frac{1}{2}} \right) + 4{\rm{L}}{{\rm{i}}_4}\left( {\frac{1}{2}} \right)\ln 2 + \frac{3}{{10}}{\ln ^5}2 + \frac{9}{4}\zeta \left( 3 \right){\ln ^2}2 - \frac{{83}}{{16}}\zeta \left( 5 \right)\nonumber\\
& \quad\quad\quad\quad\quad\quad\quad\quad- \frac{4}{3}\zeta \left( 2 \right){\ln ^3}2 - \frac{{11}}{8}\zeta \left( 2 \right)\zeta \left( 3 \right) + \frac{{11}}{4}\zeta \left( 4 \right)\ln 2.\label{c60}
\end{align}
In many other cases we are also able to obtain a formula for other Euler sum of weight five explicitly in terms of values of zeta, logarithm and polylogarithm functions. For example, we also use the method of the present paper to obtain explicit evaluations for the following sums
\begin{align}
&\sum\limits_{n = 1}^\infty  {\frac{{H_n^{\left( 3 \right)}{{\bar H}_n}}}{n}} {\left( { - 1} \right)^{n - 1}} =  - 2{\rm{L}}{{\rm{i}}_5}\left( {\frac{1}{2}} \right) + \frac{1}{{60}}{\ln ^5}2 - \frac{3}{8}\zeta \left( 3 \right){\ln ^2}2 - \frac{{193}}{{64}}\zeta \left( 5 \right) - \frac{1}{6}\zeta \left( 2 \right){\ln ^3}2 \nonumber
\\&\quad\quad\quad\quad\quad\quad\quad\quad\quad\quad+ \frac{{11}}{{16}}\zeta \left( 4 \right)\ln 2 - \frac{5}{8}\zeta \left( 2 \right)\zeta \left( 3 \right),\\
&\sum\limits_{n = 1}^\infty  {\frac{{{{\bar H}_n}\bar H_n^{\left( 3 \right)}}}{n}} {\left( { - 1} \right)^{n - 1}} = 2{\rm{L}}{{\rm{i}}_4}\left( {\frac{1}{2}} \right)\ln 2 + \frac{1}{{12}}{\ln ^5}2 + \frac{3}{8}\zeta \left( 3 \right){\ln ^2}2 - \frac{{19}}{{32}}\zeta \left( 5 \right) - \frac{1}{2}\zeta \left( 2 \right){\ln ^3}2\nonumber
\\&\quad\quad\quad\quad\quad\quad\quad\quad\quad\quad + \frac{{11}}{{16}}\zeta \left( 4 \right)\ln 2 + \frac{1}{4}\zeta \left( 2 \right)\zeta \left( 3 \right),\\
&\sum\limits_{n = 1}^\infty  {\frac{{H_n^2{{\bar H}_n}}}{{{n^2}}}}  = 4{\rm{L}}{{\rm{i}}_5}\left( {\frac{1}{2}} \right) + 2{\rm{L}}{{\rm{i}}_4}\left( {\frac{1}{2}} \right)\ln 2 + \frac{1}{{20}}{\ln ^5}2 - \frac{{165}}{{32}}\zeta \left( 5 \right) - \frac{1}{6}\zeta \left( 2 \right){\ln ^3}2 \nonumber
\\&\quad\quad\quad\quad\quad\quad+ \frac{1}{2}\zeta \left( 2 \right)\zeta \left( 3 \right) + \frac{{109}}{{16}}\zeta \left( 4 \right)\ln 2,\\
&\sum\limits_{n = {\rm{1}}}^\infty  {\frac{{H_n^2{{\bar H}_n}}}{{{n^2}}}{{\left( { - 1} \right)}^{n - 1}}}  = 2{\rm{L}}{{\rm{i}}_5}\left( {\frac{1}{2}} \right) - \frac{1}{{60}}{\ln ^5}2 - \frac{7}{8}\zeta \left( 3 \right){\ln ^2}2 - \frac{{197}}{{64}}\zeta \left( 5 \right) + \frac{1}{6}\zeta \left( 2 \right){\ln ^3}2\nonumber
\\&\quad\quad\quad\quad\quad\quad\quad\quad\quad - \frac{3}{4}\zeta \left( 2 \right)\zeta \left( 3 \right) + \frac{{109}}{{16}}\zeta \left( 4 \right)\ln 2,\\
&\sum\limits_{n = 1}^\infty  {\frac{{\bar H_n^2\bar H_n^{\left( 2 \right)}}}{n}{{\left( { - 1} \right)}^{n - 1}}}  = 14{\rm{L}}{{\rm{i}}_5}\left( {\frac{1}{2}} \right) - \frac{7}{{60}}{\ln ^5}2 - \frac{7}{8}\zeta \left( 3 \right){\ln ^2}2 - \frac{{899}}{{64}}\zeta \left( 5 \right) - \frac{2}{3}\zeta \left( 2 \right){\ln ^3}2 \nonumber
\\&\quad\quad\quad\quad\quad\quad\quad\quad\quad\quad + \frac{1}{4}\zeta \left( 2 \right)\zeta \left( 3 \right) + \frac{{193}}{{16}}\zeta \left( 4 \right)\ln 2,\\
&\sum\limits_{n = 1}^\infty  {\frac{{H_n^{\left( 2 \right)}\bar H_n^2}}{n}{{\left( { - 1} \right)}^{n - 1}}}  =  - 6{\rm{L}}{{\rm{i}}_4}\left( {\frac{1}{2}} \right)\ln 2 - \frac{1}{4}{\ln ^5}2 - \frac{1}{6}\zeta \left( 2 \right){\ln ^3}2 + \frac{1}{8}\zeta \left( 2 \right)\zeta \left( 3 \right) \nonumber
\\&\quad\quad\quad\quad\quad\quad\quad\quad\quad\quad + \frac{9}{2}\zeta \left( 4 \right)\ln 2,\\
&\sum\limits_{n = 1}^\infty  {\frac{{H_n^3{{\bar H}_n}}}{n}{{\left( { - 1} \right)}^{n - 1}}}  =  - 2{\rm{L}}{{\rm{i}}_5}\left( {\frac{1}{2}} \right) - \frac{{11}}{{60}}{\ln ^5}2 - \frac{9}{{16}}\zeta \left( 3 \right){\ln ^2}2 + \frac{{55}}{{16}}\zeta \left( 5 \right) + \frac{7}{{12}}\zeta \left( 2 \right){\ln ^3}2  \nonumber
\\&\quad\quad\quad\quad\quad\quad\quad\quad\quad\quad- \frac{1}{{16}}\zeta \left( 2 \right)\zeta \left( 3 \right) + \frac{{37}}{8}\zeta \left( 4 \right)\ln 2,\\
&\sum\limits_{n = 1}^\infty  {\frac{{{H_n}H_n^{\left( 2 \right)}{{\bar H}_n}}}{n}{{\left( { - 1} \right)}^{n - 1}}}  =  - 6{\rm{L}}{{\rm{i}}_5}\left( {\frac{1}{2}} \right) + \frac{1}{{20}}{\ln ^5}2 + \frac{7}{{16}}\zeta \left( 3 \right){\ln ^2}2 + \frac{{155}}{{32}}\zeta \left( 5 \right) - \frac{1}{{12}}\zeta \left( 2 \right){\ln ^3}2\nonumber
\\&\quad\quad\quad\quad\quad\quad\quad\quad\quad\quad\quad\quad + \frac{7}{{16}}\zeta \left( 2 \right)\zeta \left( 3 \right) - 2\zeta \left( 4 \right)\ln 2.
\end{align}
Some of the above identities can be found in \cite{Xu2016}. Moreover,
we use Mathematica tool to check numerically each of the specific identities listed. The numeric results prove that the obtained formulas are correct.

\subsection{Some evaluation of Euler-type sums}
In this subsection, we will establish some explicit relationships which
involve Euler type sums of the form
\[{S_{{m_1}{m_2}...{m_p}}}(r,k;{x_1}, \ldots ,{x_p};x) = \sum\limits_{n = 1}^\infty  {\frac{{{\zeta _n}\left( {{m_1};{x_1}} \right) \cdots {\zeta _n}\left( {{m_r};{x_p}} \right)}}{{(n + r)\left( {n + k} \right)}}}x^n,\quad p, k,r \in \N\]
and (alternating) harmonic numbers, where $x,\ x_j=\pm 1$ and $r\leq k$. We develops an approach to evaluation of Euler-type sums ${S_{{m_1}{m_2}...{m_p}}}(r,k;{x_1}, \ldots ,{x_p};x) $. The approach is based on  integral computations.

\begin{thm}\label{thm4.1} For integers $0\leq r< k$ and $i_1\leq i_2 \leq \cdots \leq i_p$, then the following relation holds:
\begin{align}\label{d1}
&\left( {k - r} \right)\left( {{S_{{i_1} \cdots {i_p}}}\left( {r,k;{x_1}, \cdots ,{x_p};1} \right) - {S_{\left( {{i_1} +  \cdots  + {i_p}} \right)}}\left( {r,k;{x_1} \cdots {x_p};1} \right)} \right)\nonumber\\
&  =  - \sum\limits_{i = 1}^{k - r} {\sum\limits_{n = 1}^\infty  {\frac{{x_1^n \cdots x_p^n}}{{{n^{{i_1} +  \cdots  + {i_p}}}\left( {n + r + i - 1} \right)}}} } \nonumber\\
&\quad - \sum\limits_{i = 1}^{k - r} {\left\{ {\sum\limits_{j = 0}^{p - 1} {\sum\limits_{1 \le {l_1} <  \cdots  < {l_j} \le p} {{{\left( { - 1} \right)}^{p - j}}\frac{{\prod\limits_{r = 1}^j {{\zeta _n}\left( {{i_{{l_r}}};{x_{{l_r}}}} \right)} }}{{{n^{{i_1} +  \cdots  + {i_p} - \left( {{i_{{l_1}}} +  \cdots  + {i_{{l_j}}}} \right)}}\left( {n + r + i - 1} \right)}}\frac{{\prod\limits_{a = 1}^p {x_a^n} }}{{\prod\limits_{b = 1}^j {x_{{i_{{l_b}}}}^n} }}} } } \right\}} .
\end{align}
\end{thm}
\pf From (\ref{c15}), we have
\begin{align}\label{d2}
&\sum\limits_{n = 1}^\infty  {\left\{ {\prod\limits_{j = 1}^p {{\zeta _n}\left( {{i_j};{x_j}} \right)}  - {\zeta _n}\left( {{i_1} +  \cdots  + {i_p};{x_1} \cdots {x_p}} \right)} \right\}{x^n}} \nonumber\\
& =  - \sum\limits_{n = 1}^\infty  {\left\{ {\frac{{\prod\limits_{j = 1}^p {x_j^n} }}{{{n^{{i_1} +  \cdots  + {i_p}}}}} + \sum\limits_{j = 0}^{p - 1} {\sum\limits_{1 \le {l_1} <  \cdots  < {l_j} \le p} {{{\left( { - 1} \right)}^{p - j}}\frac{{\prod\limits_{r = 1}^j {{\zeta _n}\left( {{i_{{l_r}}};{x_{{l_r}}}} \right)} }}{{{n^{{i_1} +  \cdots  + {i_p} - \left( {{i_{{l_1}}} +  \cdots  + {i_{{l_j}}}} \right)}}}}\frac{{\prod\limits_{a = 1}^p {x_a^n} }}{{\prod\limits_{b = 1}^j {x_{{i_{{l_b}}}}^n} }}} } } \right\}\frac{{{x^n}}}{{1 - x}}}.
\end{align}
Multiplying (\ref{d2}) by $x^{r-1}-x^{k-1}$ and integrating over (0,1), by a simple calculation we may easily deduce the desired result.\hfill$\square$\\
Hence, from Theorem \ref{thm4.1}, we can get the following corollaries.

\begin{cor} For positive integers $m_1,\ m_2,\ k$ and $r\geq 0$, we have
\begin{align}
 &\left( {k - r} \right)\sum\limits_{n = 1}^\infty  {\frac{{{H^{(m_1)} _n}{H^{(m_2)} _n} - {H^{(m_1+m_2)} _n}}}{{\left( {n + r} \right)\left( {n + k} \right)}}}\nonumber  \\
 & =\sum\limits_{j = 1}^{k - r} {\sum\limits_{n = 1}^\infty  {\left\{ {\frac{{H_n^{({m_1})}}}{{{n^{{m_2}}}\left( {n + r + j - 1} \right)}} + \frac{{H_n^{({m_2})}}}{{{n^{{m_1}}}\left( {n + r + j - 1} \right)}} - \frac{2}{{{n^{{m_1} + {m_2}}}\left( {n + r + j - 1} \right)}}} \right\}} } ,\label{d3}\\
&\left( {k - r} \right)\sum\limits_{n = 1}^\infty  {\frac{{\bar H_n^{\left( {{m_1}} \right)}\bar H_n^{\left( {{m_2}} \right)} - H_n^{\left( {{m_1} + {m_2}} \right)}}}{{\left( {n + r} \right)\left( {n + k} \right)}}} \nonumber\\
&= \sum\limits_{j = 1}^{k - r} {\sum\limits_{n = 1}^\infty  {\left\{ {\frac{{\bar H_n^{\left( {{m_1}} \right)}{{\left( { - 1} \right)}^{n - 1}}}}{{{n^{{m_2}}}\left( {n + r + j - 1} \right)}} + \frac{{\bar H_n^{\left( {{m_2}} \right)}{{\left( { - 1} \right)}^{n - 1}}}}{{{n^{{m_2}}}\left( {n + r + j - 1} \right)}} - \frac{2}{{{n^{{m_1} + {m_2}}}\left( {n + r + j - 1} \right)}}} \right\}} },\label{d4}\\
&\left( {k - r} \right)\sum\limits_{n = 1}^\infty  {\frac{{H_n^{({m_1})}\bar H_n^{\left( {{m_2}} \right)} - \bar H_n^{\left( {{m_1} + {m_2}} \right)}}}{{\left( {n + r} \right)\left( {n + k} \right)}}} \nonumber\\
& = \sum\limits_{j = 1}^{k - r} {\sum\limits_{n = 1}^\infty  {\left\{ {\frac{{H_n^{({m_1})}}{{\left( { - 1} \right)}^{n - 1}} }{{{n^{{m_2}}}\left( {n + r + j - 1} \right)}}+ \frac{{\bar H_n^{\left( {{m_2}} \right)}}}{{{n^{{m_1}}}\left( {n + r + j - 1} \right)}} - \frac{{2{{\left( { - 1} \right)}^{n - 1}}}}{{{n^{{m_1} + {m_2}}}\left( {n + r + j - 1} \right)}}} \right\}} } ,\label{d5}
\end{align}
where the Euler-type sums $\sum\limits_{n = 1}^\infty  {\frac{{{H^{(m)} _n}}}{{{n^p}\left( {n + j} \right)}}}\quad {\rm and} \sum\limits_{n = 1}^\infty  {\frac{{{{\bar H}^{(m)} _n}}}{{{n^p}\left( {n + j} \right)}}} \;\left( { p, m,j \in \N} \right)$ can be evaluated by the following formula
\[\sum\limits_{n = 1}^\infty  {\frac{{{f }\left( n \right)}}{{{n^p}\left( {n + j} \right)}}}  = \sum\limits_{i = 1}^{p - 1} {\frac{{{{\left( { - 1} \right)}^{i-1}}}}{{{j^i}}}\sum\limits_{n = 1}^\infty  {\frac{{{f }\left( n \right)}}{{{n^{p + 1 - i}}}}} }  + \frac{{{{\left( { - 1} \right)}^{p - 1}}}}{{{j^{p - 1}}}}\sum\limits_{n = 1}^\infty  {\frac{{{f }\left( n \right)}}{{n\left( {n + j} \right)}}} .\]
Note that the closed forms of the following Euler-type sums can be found in \cite{SX2017,XMZ2016,X-p2017}.
\[\sum\limits_{n = 1}^\infty  {\frac{{{H^{(m)} _n}}}{{\left( {n + r} \right)\left( {n + k} \right)}}} ,\ \sum\limits_{n = 1}^\infty  {\frac{{{{\bar H}^{(m)}_n}}}{{\left( {n + r} \right)\left( {n + k} \right)}}},\quad r\neq k .\]
\end{cor}

\begin{cor}
For integers $m_1>0,m_2>0,m_3>0$ and $k>r\ge 0$, we have
\begin{align*}
 &\left( {k - r} \right)\sum\limits_{n = 1}^\infty  {\frac{{H_n^{({m_1})}H_n^{({m_2})}H_n^{({m_3})} - H_n^{\left( {{m_1} + {m_2} + {m_3}} \right)}}}{{\left( {n + r} \right)\left( {n + k} \right)}}} \\
&  = \sum\limits_{j = 1}^{k - r} {\sum\limits_{n = 1}^\infty  {\left\{ {\frac{{H_n^{({m_1})}H_n^{({m_2})}}}{{{n^{{m_3}}}\left( {n + r + j - 1} \right)}} + \frac{{H_n^{({m_1})}H_n^{({m_3})}}}{{{n^{{m_2}}}\left( {n + r + j - 1} \right)}} + \frac{{H_n^{({m_2})}H_n^{({m_3})}}}{{{n^{{m_1}}}\left( {n + r + j - 1} \right)}}} \right\}} }  \\
 &\quad - \sum\limits_{j = 1}^{k - r} {\sum\limits_{n = 1}^\infty  {\left\{ {\frac{{H_n^{({m_1})}}}{{{n^{{m_2} + {m_2}}}\left( {n + r + j - 1} \right)}} + \frac{{H_n^{({m_2})}}}{{{n^{{m_1} + {m_3}}}\left( {n + r + j - 1} \right)}} + \frac{{H_n^{({m_3})}}}{{{n^{{m_1} + {m_2}}}\left( {n + r + j - 1} \right)}}} \right\}} } .
 \end{align*}
\end{cor}

It is possible that some of other Euler type sums can be obtained using techniques of the present paper.

{\bf Acknowledgments.} The authors would like to thank the anonymous
referee for his/her helpful comments, which improve the presentation
of the paper.
 \end{document}